\documentclass{article}
\usepackage{amsfonts}

\usepackage[dvips]{graphicx}
\usepackage{amsmath}

\newtheorem{theorem}{Theorem}[subsection]

\newtheorem{conjecture}[theorem]{Conjecture}
\newtheorem{corollary}[theorem]{Corollary}

\newtheorem{definition}[theorem]{Definition}
\newtheorem{example}[theorem]{Example}

\newtheorem{lemma}[theorem]{Lemma}

\newtheorem{proposition}[theorem]{Proposition}

\newenvironment{proof}[1][Proof]{\textbf{#1.} }{\ \rule{0.5em}{0.5em}}

\input{tcilatex}

\begin{document}

\title{Quantization of branching coefficients for classical Lie groups}
\author{C\'{e}dric Lecouvey \\
Laboratoire de Math\'{e}matiques Pures et Appliqu\'{e}es Joseph Liouville\\
Universit\'{e} du Littoral C\^{o}te d'Opale\\
Maison de la recherche Blaise Pascal\\
50 rue F. Buisson B.P. 699 62228 Calais Cedex}
\date{}
\maketitle

\begin{abstract}
We study natural quantizations of branching coefficients corresponding to
the restrictions of the classical Lie groups to their Levi subgroups.\ We
show that they admit a stable limit which can be regarded as a $q$-analogue
of a tensor product multiplicity. According to a conjecture by Shimozono,
the stable one-dimensional sum for nonexceptional affine crystals are
expected to occur as special cases of these $q$-analogues.
\end{abstract}

\section{Introduction}

The Kostka coefficients and the Littlewood-Richardson coefficients which
have many occurrences in the representation theory of $GL_{n}$ admit
interesting $q$-analogues. Giving $\lambda $ a partition of length at most $%
n $ and $\mu \in \mathbb{Z}^{n}$, the $q$-analogue of the Kostka coefficient 
$K_{\lambda ,\mu }$ giving the dimension of the weight space $\mu $ in the
irreducible finite dimensional $GL_{n}$-module $V^{GL_{n}}(\lambda )$ of
highest weight $\lambda $ is the Kostka-Foulkes polynomial $K_{\lambda ,\mu
}(q)$ (also called Lusztig $q$-analogue of weight multiplicity). Consider $%
\mathbf{\mu }=(\mu ^{(1)},...,\mu ^{(r)})$ a $r$-tuple of partitions of
lengths summing $n$ and denote by $\mu $ the $n$-tuple obtained by reading
the parts of the $\mu ^{(p)}$'s from left to right.\ There exist in the
literature different quantizations of the Littlewood-Richardson coefficient $%
c_{\mu ^{(1)},...,\mu ^{(r)}}^{\lambda }$ giving the multiplicity of $%
V^{GL_{n}}(\lambda )$ in the tensor product $V^{GL_{n}}(\mu ^{(1)})\otimes
\cdot \cdot \cdot \otimes V^{GL_{n}}(\mu ^{(r)}).$

\noindent In \cite{LLT0} Lascoux Leclerc and Thibon have introduced such a $%
q $-analogue by mean of certain generalizations of semi-standard Young
tableaux called ribbon tableaux. They have proved in \cite{LLT} that the
polynomials obtained belong to a family of parabolic Kazhdan-Lusztig
polynomials introduced by Deodhar which have nonnegative integer
coefficients \cite{KAT}.

\noindent When $\mu ^{(1)},...,\mu ^{(r)}$ are rectangular partitions, it is
also possible to define $q$-analogues of the coefficients $c_{\mu
^{(1)},...,\mu ^{(r)}}^{\lambda }$ by considering the one-dimensional sums $%
X_{\lambda ,\mathbf{\mu }}^{\emptyset }(q)$ obtained from affine $%
A_{n}^{(1)} $-crystals associated to Kirillov-Reshetikhin $U_{q}^{\prime }(%
\widehat{sl_{n}})$-modules \cite{HKOTY}.

\noindent Consider $\eta =(\eta _{1},...,\eta _{r})$ a sequence of positive
integers summing $n$ and suppose that $\mu ^{(p)}$ has length $\eta _{p}$
for any $p=1,...,r.\;$The Littlewood-Richardson coefficient $c_{\mu
^{(1)},...,\mu ^{(r)}}^{\lambda }$ also coincide with the multiplicity of
the tensor product $V^{GL_{\eta _{1}}}(\mu ^{(1)})\otimes \cdot \cdot \cdot
\otimes V^{GL_{\eta _{r}n}}(\mu ^{(r)})$ in the restriction of $%
V^{GL_{n}}(\lambda )$ to its Levi subgroup $GL_{\eta }=GL_{\eta _{1}}\times
\cdot \cdot \cdot \times GL_{\eta _{r}}.$ This duality permits to express $%
c_{\mu ^{(1)},...,\mu ^{(r)}}^{\lambda }$ in terms of a Kostant-type
partition function.\ By quantifying this partition function, Shimozono and
Weyman \cite{SZ} have introduced another natural $q$-analogue of $c_{\mu
^{(1)},...,\mu ^{(r)}}^{\lambda }$ that we will denote $K_{\lambda ,\mu
}^{GL_{n},I}(q)$ ($I$ being the set of the simple roots of $GL_{\eta }$).\
The polynomials $K_{\lambda ,\mu }^{GL_{n},I}(q)$ are Poincar\'{e}
polynomials and appear in the Hilbert series of the Euler characteristic of
certain graded virtual $G$-modules. By a result of Broer \cite{broer}, they
admit nonnegative coefficients providing that the $\mu ^{(p)}$'s are
rectangular partitions of decreasing heights.\ In this case Shimozono has
proved in \cite{Sh0} that $K_{\lambda ,\mu }^{GL_{n},I}(q)$ coincide with $%
X_{\lambda ,\mathbf{\mu }}^{\emptyset }(q)$. This result which is based on a
combinatorial description of the polynomials $K_{\lambda ,\mu
}^{GL_{n},I}(q),$ permits in particular to recover that they have
nonnegative coefficients independently of the results of Broer.\ Under the
same hypothesis, it is conjectured that $K_{\lambda ,\mu }^{GL_{n},I}(q)$
also coincide with the LLT quantization of $c_{\mu ^{(1)},...,\mu
^{(r)}}^{\lambda }.$ When the $\mu ^{(p)}$'s are simply row partitions, we
have $\mu =(\mu ^{(1)},...,\mu ^{(r)})\in \mathbb{Z}^{n}$ and $K_{\lambda
,\mu }^{GL_{n},I}(q)$ is the Kostka-Foulkes polynomial associated to the
weights $\lambda $ and $\mu $.

\bigskip

\noindent Let $G$ be one of the classical groups $GL_{n},SO_{2n+1},Sp_{2n}$
or $SO_{2n}$ and $R_{G}^{+}$\ its set of positive roots.\ Kostka and
Littlewood-Richardson coefficients can be regarded as branching coefficients
corresponding to the restriction of $GL_{n}$ to its principal Levi
subgroups.\ This naturally yields to study the branching coefficients
corresponding to the restriction to a subgroup $G_{0}$ (not necessarily of
Levi type) of $G$ and their corresponding $q$-analogues.\ The branching
coefficients which are considered in this paper can be expressed in terms of
certain partition functions counting the number of way to decompose a weight
of $G$ into a linear positive combination of simple roots belonging to a
fixed subset of $R_{G}^{+}$.\ To obtain the corresponding $q$-analogues of
these coefficients, it suffices to quantify these partitions functions. In
particular when $G_{0}=H_{G}$ is the maximal torus of $G,$ the $q$-analogues
obtained in this way are precisely the Lusztig $q$-analogues of weight
multiplicities associated to $G$ \cite{Lu}.

\noindent Our aim in this paper is two fold.\ First, we study the natural
quantizations of branching coefficients corresponding to the restrictions of 
$G$ to the Levi subgroups of its standard parabolic subgroups.\ These
polynomials will be denoted $K_{\lambda ,\mu }^{G,I}(q)$ where $I$ is the
set of simple roots of the Levi subgroup $L_{G,I}$ considered.\ The
polynomials $K_{\lambda ,\mu }^{G,I}(q)$ are generalizations of Lusztig $q$%
-analogues of weight multiplicities which coincide with the $q$-analogues $%
K_{\lambda ,\mu }^{GL_{n},I}(q)$ for $G=GL_{n}$. From the results of \cite
{broer}, one can derive that they have nonnegative coefficients when $\mu $
is stable under the action of the Weyl group of $L_{G,I}.\;$The polynomials $%
K_{\lambda ,\mu }^{GL_{n},I}(q)$ indexed by pairs of dominant weights $%
\lambda $ and $\mu $ which contain sufficiently large multiples of the
fundamental weight $\kappa =(1,...,1)\in \mathbb{N}^{n}$ are called stable.\
This terminology reflects the fact that they are invariant if $\lambda $ and 
$\mu $ are translated by $\kappa $.\ When the Levi subgroup $L_{G,I}$ is
isomorphic to a direct product of linear groups, we prove that this stable
limit decomposes as nonnegative integer combination of polynomials $%
K_{\lambda ,\mu }^{GL_{n},I}(q)$ (Theorem \ref{prop_dec_K_c}). This result
can be regarded as a generalization of the decomposition of the stable limit
of Lusztig $q$-analogues associated to $G$ as a sum of Kostka-Foulkes
polynomials given in \cite{lec}.\ For a general Levi subgroup, we conjecture
that the polynomials $K_{\lambda ,\mu }^{G,I}(q)$ have nonnegative
coefficients providing that $\mu $ is a partition (Conjecture \ref{conj}).
Note that this condition is in particular fulfilled when $\mu $ is stable
under the action of the Weyl group of $L_{G,I}.$ While writing this paper,
the author was informed that an equivalent statement of this conjecture
first appeared in some unpublished notes by Broer.

\noindent Next we study Littlewood-Richardson-type coefficients associated
to $G=SO_{2n+1},Sp_{2n}$ or $SO_{2n}$ and discuss the problem of their
possible quantizations.$\;$Note first that the LLT quantization of the
Littlewood-Richardson coefficients for\ $GL_{n}$ is based on a very special
property of the plethysm of the Schur functions with the power sums.\ Indeed
the coefficients appearing in the decomposition of this plethysm on the
basis of the Schur functions are, up to a sign, Littlewood-Richardson
coefficients. An analogous property for the other classical groups does not
exist.\ Thus it seems impossible to relate $q$-analogues of tensor product
multiplicities to Deodhar's polynomials by proceeding as in \cite{LLT}. With
the above notation for $\mathbf{\mu }$ and $\eta $, we define the
coefficient $d_{\mu ^{(1)},...,\mu ^{(r)}}^{\lambda }$ as the multiplicity
of the finite dimensional irreducible $G$-module $V^{G}(\lambda )$ of
highest weight $\lambda $ in the tensor product $V^{G}(\mu ^{(1)})\otimes
\cdot \cdot \cdot \otimes V^{G}(\mu ^{(r)}).$ We show that this coefficient
can be expressed in terms of a partition function (Proposition \ref{pro_d}%
).\ This implies in particular that it does not depend on $G$.\ We also
establish a duality result (Proposition \ref{th_dual1}) between the
coefficients $d_{\mu ^{(1)},...,\mu ^{(r)}}^{\lambda }$ and certain
branching coefficients corresponding to the restriction of $SO_{2n}$ to the
subgroup $SO_{2\eta _{1}}\times \cdot \cdot \cdot \times SO_{2\eta _{r}}$
(which is not isomorphic to a Levi subgroup of $SO_{2n}$). This permits to
define $q$-analogues for the coefficients $d_{\mu ^{(1)},...,\mu
^{(r)}}^{\lambda }$ but the polynomials obtained in this way may have
negative coefficients even if the $n$-tuple $\mu $ is a partition. Denote by 
$\frak{V}^{G}(\lambda )$ the restriction of the irreducible finite
dimensional $GL_{N}$-module of highest weight $\lambda $ to $G$ where $%
N=2n+1 $ if $G=SO_{2n+1}$ and $N=2n$ if $G=Sp_{2n}$ or $SO_{2n}.$ By
replacing each module $V^{G}(\mu ^{(p)})$ by $\frak{V}^{G}(\mu ^{(p)})$ in
the definition of $d_{\mu ^{(1)},...,\mu ^{(r)}}^{\lambda }$, one obtains
tensor product coefficients $\frak{D}_{\mu ^{(1)},...,\mu ^{(r)}}^{\lambda
,G}$ which can also be expressed in terms of a partition function.\ Thus,
they admit natural quantizations $\frak{D}_{\mu ^{(1)},...,\mu
^{(r)}}^{\lambda ,G}(q)$.\ Note that this time the coefficients $\frak{D}%
_{\mu ^{(1)},...,\mu ^{(r)}}^{\lambda ,G}$ and the polynomials $\frak{D}%
_{\mu ^{(1)},...,\mu ^{(r)}}^{\lambda ,G}(q)$ depend on the lie group $G$
considered.\ We obtain a duality between the $q$-analogues $\frak{D}_{\mu
^{(1)},...,\mu ^{(r)}}^{\lambda ,G}(q)$ and the stable limit of the
polynomials $K_{\lambda ,\mu }^{G,I}(q)$ associated to the Levi subgroup $%
GL_{\eta _{1}}\times \cdot \cdot \cdot \times GL_{\eta _{r}}$. In particular
the polynomials $\frak{D}_{\mu ^{(1)},...,\mu ^{(r)}}^{\lambda ,G}(q)$
decomposes as nonnegative integer combination of polynomials $K_{\lambda
,\mu }^{GL_{n},I}(q)$ (Theorem \ref{th_qdual}) and have nonnegative integer
coefficients when the $\mu ^{(p)}$'s are rectangular partitions of
decreasing heights. Within each nonexceptional family of affine algebras,
the one-dimensional sums have large rank limits which are called stable
one-dimensional sums \cite{SZa}. There exist four distinct kinds of stable
one-dimensional sums $X^{\diamondsuit }$ labelled by the symbols $%
\diamondsuit =\emptyset ,(1),(2),(1,1)$.\ The stable one-dimensional sums of
kind $\emptyset $ are related to $A_{n-1}^{(1)}$-affine crystals whereas the
stable one-dimensional sums of kind $(1),(2),(1,1)$ are defined from the
other nonexceptional families of affine crystals. Then, according to
Conjecture $5$ of \cite{sh} giving the decomposition of $X^{\diamondsuit }$
in terms of one-dimensional sums $X^{\emptyset }$, the three families of $q$%
-analogues $\frak{D}_{\mu ^{(1)},...,\mu ^{(r)}}^{\lambda ,G}(q),$ $%
G=SO_{2n+1},Sp_{2n}$ and $SO_{2n}$ should coincide (up to a simple
renormalization) respectively with the stable one-dimensional sums of kind $%
(1),(1,1)$ and $(2)$ associated to $\mu $ when the $\mu ^{(p)}$'s are
rectangular partitions of decreasing heights. This means that it should be
possible to extend the results of \cite{LS} which holds when the $\mu ^{(p)}$%
's are row partitions of decreasing heights to all stable one-dimensional
sums by establishing Conjecture 5 of \cite{sh}.

\bigskip

The paper is organized as follows.\ In Section $2$ we review the necessary
background on branching multiplicity formulas and Levi-subgroups for
classical Lie groups.\ In particular we introduce the Kostant-type partition
functions which permit to compute the branching coefficients we use in the
sequel. Section $3$ is concerned with the $q$-analogues of branching
coefficients corresponding to the restrictions to Levi subgroups.\ We prove
that they admit a stable limit which decompose as nonnegative integer
combination of Poincar\'{e} polynomials when the Levi subgroup considered is
isomorphic to a direct product of linear groups. In section $4,$ we use the
Jacobi-Trudi type determinantal expressions for the Schur functions of
classical type to derive a duality between the Littlewood-Richardson
coefficients $d_{\mu ^{(1)},...,\mu ^{(r)}}^{\lambda }$ and the branching
coefficients corresponding to the restriction of $SO_{2n}$ to $SO_{2\eta
_{1}}\times \cdot \cdot \cdot \times SO_{2\eta _{r}}.$ This duality and the
arguments used to proved it generalize the results of \cite{lec}
(corresponding to the case when all the $\mu ^{(p)}$'s are row partitions).
We observe that the natural quantization of the multiplicities $d_{\mu
^{(1)},...,\mu ^{(r)}}^{\lambda }$ may have negative coefficients. We then
introduce the polynomials $\frak{D}_{\mu ^{(1)},...,\mu ^{(r)}}^{\lambda
,G}(q)$ and show how they are related to the $q$-analogues of the branching
coefficients corresponding to the restriction of $G$ to $GL_{\eta
_{1}}\times \cdot \cdot \cdot \times GL_{\eta _{r}}.$

\section{Background}

\subsection{Branching multiplicity formulas \label{deG0}}

In the sequel $G$ is one of the complex Lie groups $GL_{n},Sp_{2n},SO_{2n+1}$
or $SO_{2n}$ and $\frak{g}$ its Lie algebra.\ We follow the convention of 
\cite{KT} to realize $G$ as a subgroup of $GL_{N}$ and $\frak{g}$\ as a
subalgebra of $\frak{gl}_{N}$ where 
\begin{equation*}
N=\left\{ 
\begin{tabular}{l}
$n$ when $G=GL_{n}$ \\ 
$2n$ when $G=Sp_{2n}$ \\ 
$2n+1$ when $G=SO_{2n+1}$ \\ 
$2n$ when $G=SO_{2n}$%
\end{tabular}
\right. .
\end{equation*}
With this convention the maximal torus $T_{G}$ of $G$ and the Cartan
subalgebra $\frak{h}_{G}$ of $\frak{g}$ coincide respectively with the
subgroup and the subalgebra of diagonal matrices of $G$ and $\frak{g}$.
Similarly the Borel subgroup $B_{G}$ of $G$ and the Borel subalgebra $\frak{b%
}_{G}$ of $\frak{g}$ coincide respectively with the subgroup and subalgebra
of upper triangular matrices of $G$ and $\frak{g}$.

\noindent Let $d_{N}$ be the linear subspace of $\frak{gl}_{N}$ consisting
of the diagonal matrices.\ For any $i\in \{1,...,n\},$ write $\varepsilon
_{i}$ for the linear map $\varepsilon _{i}:d_{N}\rightarrow \mathbb{C}$ such
that $\varepsilon _{i}(D)=\delta _{i}$ for any diagonal matrix $D$ whose $%
(i,i)$-coefficient is $\delta _{i}.$ Then $(\varepsilon _{1},...,\varepsilon
_{n})$ is an orthonormal basis of the Euclidean space $\frak{h}_{G,\mathbb{R}%
}^{\ast }$ (the real part of $\frak{h}_{G}^{\ast }).$ Let $R_{G}$ be the
root system associated to $G.\;$We can take for the simple roots of $\frak{g}
$%
\begin{equation}
\left\{ 
\begin{tabular}{l}
$\Sigma _{GL(n)}^{+}=\{\alpha _{i}=\varepsilon _{i}-\varepsilon _{i+1}\text{%
, }i=1,...,n-1\text{ for the root system }A_{n-1}\}$ \\ 
$\Sigma _{SO_{2n+1}}^{+}=\{\alpha _{n}=\varepsilon _{n}\text{ and }\alpha
_{i}=\varepsilon _{i}-\varepsilon _{i+1}\text{, }i=1,...,n-1\text{ for the
root system }B_{n}\}$ \\ 
$\Sigma _{Sp_{2n}}^{+}=\{\alpha _{n}=2\varepsilon _{n}\text{ and }\alpha
_{i}=\varepsilon _{i}-\varepsilon _{i+1}\text{, }i=1,...,n-1\text{ for the
root system }C_{n}\}$ \\ 
$\Sigma _{SO_{2n}}^{+}=\{\alpha _{n}=\varepsilon _{n}+\varepsilon _{n-1}%
\text{ and }\alpha _{i}=\varepsilon _{i}-\varepsilon _{i+1}\text{, }%
i=1,...,n-1\text{ for the root system }D_{n}\}$%
\end{tabular}
\right. .  \label{simple_roots}
\end{equation}
Then the set of positive roots are 
\begin{equation*}
\left\{ 
\begin{tabular}{l}
$R_{GL_{n}}^{+}=\{\varepsilon _{i}-\varepsilon _{j},\varepsilon
_{i}+\varepsilon _{j}\text{ with }1\leq i<j\leq n\}\text{ for the root
system }A_{n-1}$ \\ 
$R_{SO_{2n+1}}^{+}=\{\varepsilon _{i}-\varepsilon _{j},\varepsilon
_{i}+\varepsilon _{j}\text{ with }1\leq i<j\leq n\}\cup \{\varepsilon _{i}%
\text{ with }1\leq i\leq n\}\text{ for the root system }B_{n}$ \\ 
$R_{Sp_{2n}}^{+}=\{\varepsilon _{i}-\varepsilon _{j},\varepsilon
_{i}+\varepsilon _{j}\text{ with }1\leq i<j\leq n\}\cup \{2\varepsilon _{i}%
\text{ with }1\leq i\leq n\}\text{ for the root system }C_{n}$ \\ 
$R_{SO_{2n}}^{+}=\{\varepsilon _{i}-\varepsilon _{j},\varepsilon
_{i}+\varepsilon _{j}\text{ with }1\leq i<j\leq n\}\text{ for the root
system }D_{n}$%
\end{tabular}
\right. .
\end{equation*}
We denote by $R_{G}$ the set of roots of $G.\;$The Weyl group of $GL_{n}$ is
the symmetric group $S_{n}$ and for $G=SO_{2n+1},Sp_{2n}$ or $SO_{2n},$ the
Weyl group $W_{G}$ of the Lie group $G$ is the subgroup of the permutation
group of the set $\{\overline{n},...,\overline{2},\overline{1},1,2,...,n\}$\
generated by the permutations 
\begin{equation*}
\left\{ 
\begin{tabular}{l}
$s_{i}=(i,i+1)(\overline{i},\overline{i+1}),$ $i=1,...,n-1$ and $s_{n}=(n,%
\overline{n})$ $\text{for the root systems }B_{n}$ and $C_{n}$ \\ 
$s_{i}=(i,i+1)(\overline{i},\overline{i+1}),$ $i=1,...,n-1$ and $%
s_{n}^{\prime }=(n,\overline{n-1})(n-1,\overline{n})$ $\text{for the root
system }D_{n}$%
\end{tabular}
\right.
\end{equation*}
where for $a\neq b$ $(a,b)$ is the simple transposition which switches $a$
and $b.$ For $G=SO_{2n+1},Sp_{2n}$ or $SO_{2n},$ we identify the subgroup of 
$W_{G}$ generated by $s_{i}=(i,i+1)(\overline{i},\overline{i+1}),$ $%
i=1,...,n-1$ with the symmetric group $S_{n}.$ We denote by $\ell $ the
length function corresponding to the above set of generators. The action of $%
w\in W_{G}$ on $\beta =(\beta _{1},...,\beta _{n})\in \frak{h}_{G,\mathbb{R}%
}^{\ast }$ is defined by 
\begin{equation*}
w\cdot (\beta _{1},...,\beta _{n})=(\beta _{1}^{w},...,\beta _{n}^{w})
\end{equation*}
where $\beta _{i}^{w}=\beta _{w(i)}$ if $\sigma (i)\in \{1,...,n\}$ and $%
\beta _{i}^{w}=-\beta _{w(\overline{i})}$ otherwise. We denote by $\rho _{G}$
the half sum of the positive roots of $R_{G}^{+}$.\ The dot action of $W_{G}$
on $\beta =(\beta _{1},...,\beta _{n})\in \frak{h}_{\mathbb{R}}^{\ast }$ is
defined by 
\begin{equation}
w\circ \beta =w\cdot (\beta +\rho _{G})-\rho _{G}.  \label{dotaction}
\end{equation}
Write $P_{G}^{+}$ for the cone of dominant weights of $G.\;$Denote by $%
\mathcal{P}_{n}$ the set of partitions with at most $n$ parts. Each
partition $\lambda =(\lambda _{1},...,\lambda _{n})\in $ $\mathcal{P}_{n}$
will be identified with the dominant weight $\sum_{i=1}^{n}\lambda
_{i}\varepsilon _{i}.\;$Then the irreducible finite dimensional
representations of $G$ are parametrized by the partitions of $\mathcal{P}%
_{n} $.\ For any $\lambda \in \mathcal{P}_{n},$ denote by $V^{G}(\lambda )$
the irreducible finite dimensional representation of $G$ of highest weight $%
\lambda .$ In the sequel we will also need the irreducible rational
representations of $GL_{n}$.\ They are indexed by the $n$-tuples 
\begin{equation}
(\gamma ^{+},\gamma ^{-})=(\gamma _{1}^{+},\gamma _{2}^{+},...,\gamma
_{p}^{+},0,...,0,-\gamma _{q}^{-},...,-\gamma _{1}^{-})  \label{jamma+-}
\end{equation}
where $\gamma ^{+}$ and $\gamma ^{-}$ are partitions of length $p$ and $q$
such that $p+q\leq n.$ Write $\widetilde{\mathcal{P}}_{n}$ for the set of
such $n$-tuples and denote also by $V^{GL_{n}}(\gamma )$ the irreducible
rational representation of $GL_{n}$ of highest weight $\gamma =(\gamma
^{+},\gamma ^{-})\in \widetilde{\mathcal{P}}_{n}.$

\noindent Consider a $r$-tuple $\eta =(\eta _{1},...,\eta _{r})$ of positive
integers summing $n.$ Given $\mathbf{\mu }=(\mu ^{(1)},...,\mu ^{(r)})$ a $r$%
-tuple of partitions such that $\mu ^{(p)}\in \mathcal{P}_{\eta _{p}}$ for $%
p=1,...,r,$ we denote by $\mu $ the $n$-tuple obtained by reading
successively the parts of the partitions $\mu ^{(1)},...,\mu ^{(r)}$ from
left to right.

\bigskip

\noindent As customary, we identify $P_{G}$ the lattice of weights of $G$
with a sublattice of $(\frac{1}{2}\mathbb{Z})^{n}.$ For any $\beta =(\beta
_{1},...,\beta _{n})\in P_{G},$ we set $\left| \beta \right| =\beta
_{1}+\cdot \cdot \cdot +\beta _{n}.\;$We use for a basis of the group
algebra $\mathbb{Z}[\mathbb{Z}^{n}],$ the formal exponentials $(e^{\beta
})_{\beta \in \mathbb{Z}^{n}}$ satisfying the relations $e^{\beta
_{1}}e^{\beta _{2}}=e^{\beta _{1}+\beta _{2}}.$ We furthermore introduce $n$
independent indeterminates $x_{1},...,x_{n}$ in order to identify $\mathbb{Z}%
[\mathbb{Z}^{n}]$ with the ring of polynomials $\mathbb{Z}%
[x_{1},...,x_{n},x_{1}^{-1},...,x_{n}^{-1}]$ by writing $e^{\beta
}=x_{1}^{\beta _{1}}\cdot \cdot \cdot x_{n}^{\beta _{n}}=x^{\beta }$ for any 
$\beta =(\beta _{1},...,\beta _{n})\in \mathbb{Z}^{n}.$

\noindent For any $\lambda \in \mathcal{P}_{n},$ we denote by $s_{\lambda
}^{G}$ the universal character of type $G$ associated to $\lambda $ and by $%
\mathcal{F}^{G}$ the ring of the universal characters of type $G$ defined by
Koike and Terada \cite{KAT}.

\bigskip

\noindent We now give a technical lemma that we will be led to use in the
sequel.\ Consider $\lambda \in \mathcal{P}_{n}$, $\mu \in \mathbb{Z}^{n}$
and $G$ one of the Lie groups $SO_{2n+1},Sp_{2n}$ or $SO_{2n}.$ Set $\kappa
=(1,...,1)\in \mathbb{N}^{n}.$

\begin{lemma}
\label{lem_tech}Let $\mathcal{M}:\mathbb{Z}^{n}\rightarrow \mathbb{Z}[q]$ be
a map such that for any $\beta \in \mathbb{Z}^{n}$, $\mathcal{M}(\beta )=0$
if $\left| \beta \right| <0.$ Then for any integer $k\geq \frac{\left|
\lambda \right| -\left| \mu \right| }{2}$ we have: 
\begin{equation}
\sum_{w\in W_{G}}(-1)^{\ell (w)}\mathcal{M}(w(\lambda +k\kappa +\rho
_{G})-\mu -k\kappa -\rho _{G})=\sum_{\sigma \in S_{n}}(-1)^{\ell (\sigma )}%
\mathcal{M}(\sigma (\lambda +\rho )-\mu -\rho )  \label{eq}
\end{equation}
where $\rho =(n,n-1,...,1).$
\end{lemma}

\begin{proof}
Consider $\delta =(\delta _{1},...,\delta _{n})\in \mathbb{Z}^{n}$ and $w\in
W_{G}.\;$Write $w(\delta )=(\delta _{1}^{w},...,\delta _{n}^{w})$ and set $%
E_{w}=\{i\mid w(i)\notin \{1,...,n\}\}$. Define the sum $S_{w,\delta
}=\sum_{i\in E_{w}}\delta _{i_{k}}.\;$Then $\left| w(\delta )\right| =\left|
\delta \right| -2S_{w,\delta }$.\ Now consider $k$ a nonnegative integer and
set $\delta =(\lambda +\rho _{G}+k\kappa ).\;$We have $\left| w(\lambda
+\rho _{G}+k\kappa )\right| =\left| (\lambda +\rho _{G}+k\kappa )\right|
-2S_{w,\delta }.\;$But $S_{w,\delta }=S_{w,\lambda +\rho _{G}}+kp$ where $p=%
\mathrm{card}(E_{w}).\;$Thus we obtain 
\begin{multline*}
\left| w(\lambda +\rho _{G}+k\kappa )-(\mu +\rho _{G}+k\kappa )\right|
=\left| (\lambda +\rho _{G}+k\kappa )\right| -2S_{w,\lambda +\rho
_{G}}-\left| (\mu +\rho _{G}+k\kappa )\right| -2kp= \\
\left| \lambda \right| -\left| \mu \right| -2S_{w,\lambda +\rho _{G}}-2kp.
\end{multline*}
When $w\notin \mathcal{S}_{n},$ we have $p\geq 1\;$and $S_{w,\lambda +\rho
_{G}}\geq 1$ since the coordinates of $\lambda +\rho _{G}$ are all positive.
Hence $\left| w(\lambda +\rho _{G}+k\kappa )-(\mu +\rho _{G}+k\kappa
)\right| <\left| \lambda \right| -\left| \mu \right| -2k$ and is negative as
soon as $k\geq \frac{\left| \lambda \right| -\left| \mu \right| }{2}.$ For
such an integer $k$ the sum defining the left hand side of the equality (\ref
{eq}) normally running over $W_{G}$ can be restricted to $S_{n}$.\ Moreover
we can write $\rho _{G}=\rho +\varepsilon \kappa $ with $\varepsilon =-\frac{%
1}{2},0$ or $1$ respectively for $G=SO_{2n+1},Sp_{2n}$ or $SO_{2n}.$ Since $%
\sigma (p\kappa )=p\kappa $ for any $p\in (\frac{1}{2}\mathbb{Z}$)$^{n},$
this yields to the desired equality.
\end{proof}

\bigskip

\noindent Let $G_{0}\subset G$ be a complex Lie subgroup of $G$ and $\frak{g}%
_{0}$ its Lie algebra. We suppose in the sequel that $G_{0}$ is isomorphic
to a product of classical Lie groups whose maximal torus $T_{0}$ is equal to 
$T_{G}.$ Let $\frak{h}_{0}$ be the Cartan subalgebra of $\frak{g}_{0}.$ We
have $\frak{h}_{0}=\frak{h}_{G}\frak{,}$ thus $\frak{h}_{0}^{\ast }=\frak{h}%
_{G}^{\ast }$.\ In particular we can consider the set $R_{0}$ of roots of $%
\frak{g}_{0}$ and its subset of positive roots $R_{0}^{+}$ respectively as
subsets of $R_{G}$ and $R_{G}^{+}.$

\noindent The partition function $\mathcal{P}_{G_{0}}$ associated to $G_{0}$
is defined by the formal identity 
\begin{equation}
\prod_{\alpha \in R^{+}-R_{0}^{+}}\frac{1}{1-e^{\alpha }}=\sum_{\beta \in
P_{G}}\mathcal{P}^{G_{0}}(\beta )e^{\beta }.  \label{def_part}
\end{equation}
Note that $\mathcal{P}_{G_{0}}$ coincide with the Kostant partition function
when $G_{0}=T_{G}$ (that is $R_{0}^{+}=\emptyset $)$.$ Write $P_{G_{0}}^{+}$
for the cone of dominant weights of $G_{0}.$ For any $\lambda $ in $%
P_{G}^{+} $ and $\mathbf{\mu }$ in $P_{G_{0}}^{+}$ we denote by $[V(\lambda
)^{G}:V(\mathbf{\mu })^{G_{0}}]$ the multiplicity of the irreducible $G_{0}$%
-module $V(\mathbf{\mu })^{G_{0}}$ of highest weight $\mathbf{\mu }$ in the
restriction of the $G$-module $V(\lambda )^{G}$ to $G_{0}.$

\begin{theorem}
\label{th_branch}With the above notation we have 
\begin{equation*}
\lbrack V(\lambda )^{G}:V(\mathbf{\mu })^{G_{0}}]=\sum_{w\in
W_{G}}(-1)^{\ell (w)}\mathcal{P}^{G_{0}}(w\circ \lambda -\mathbf{\mu }).
\end{equation*}
\end{theorem}

\noindent We refer the reader to Theorem 8.2.1 of \cite{GW} for the proof.

\subsection{Branching coefficients associated to Levi subgroups\label%
{subsec_levi}}

Consider $I$ a subset of $\Sigma _{G}^{+}$ the set of simple roots
associated to the classical Lie algebra $G.$ Denote by $\pi _{G,I}$ the
standard parabolic subgroup of $G$ (that is containing the Borel subgroup $%
B_{G})$ defined by $I.$ Recall that the roots of $\pi _{G,I}$ are those of $%
R_{G}^{+}$ together with the negative roots of $R_{G}$ which are $\mathbb{Z}$%
-linear combinations of the simple roots contained in $I.\;$Write $L_{G,I}$
for the Levi subgroup of the parabolic $\pi _{G,I}$ and $\frak{l}_{G,I}$ its
corresponding Lie algebra.\ Let $R_{G,I}$ be the subsystem of roots spanned
by $I$ and $R_{G,I}^{+}$ the subset of positive roots in $R_{G,I}.$ Then $%
R_{G,I}$ and $R_{G,I}^{+}$ are respectively the set of roots and the set of
positive roots of $\frak{l}_{G,I}$.

\noindent The Levi subgroup $L_{G,I}$ corresponds to the removal, in the
Dynkin diagram of $G,$ of the nodes which are not associated to a simple
root belonging to $I$. Write 
\begin{equation*}
J=\Sigma _{G}^{+}-I=\{\alpha _{j_{1}},...,\alpha _{j_{r}}\}
\end{equation*}
where for any $k=1,...,r,$ $\alpha _{j_{k}}$ is a simple root of $\Sigma
_{G}^{+}$ and $j_{1}<\cdot \cdot \cdot <j_{r}.$ Set $l_{1}=j_{1},$ $%
l_{k}=j_{k}-j_{k-1},$ $k=2,...,r$ and $l_{r+1}=n-j_{r}.$ According to \cite
{K}, the Levi group $L_{G,I}$ is isomorphic to a direct product of classical
Lie groups determined by the $(r+1)$-tuple $l_{I}=(l_{1},...,l_{r+1})$ of
nonnegative integers summing $n.$ We give in the table below the direct
product associated to each Levi group $L_{G,I}.$%
\begin{equation}
\begin{tabular}{|l|l|}
\hline
$1:G=GL_{n}$ & $L_{G,I}\simeq GL_{l_{1}}\times \cdot \cdot \cdot \times
GL_{l_{r+1}}$ \\ 
$2:G=SO_{2n+1}$ and $l_{r+1}\geq 2$ & $L_{G,I}\simeq GL_{l_{1}}\times \cdot
\cdot \cdot \times GL_{l_{r}}\times SO_{2l_{r+1}+1}$ \\ 
$3:G=Sp_{2n}$ and $l_{r+1}\geq 2$ & $L_{G,I}\simeq GL_{l_{1}}\times \cdot
\cdot \cdot \times GL_{l_{r}}\times Sp_{2l_{r+1}}$ \\ 
$4:G=SO_{2n}$ and $l_{r+1}\geq 4$ & $L_{G,I}\simeq GL_{l_{1}}\times \cdot
\cdot \cdot \times GL_{l_{r}}\times SO_{2l_{r+1}}$ \\ 
$5:G=SO_{2n+1},Sp_{2n}$ or $SO_{2n}$ and $l_{r+1}=0$ & $L_{G,I}\simeq
GL_{l_{1}}\times \cdot \cdot \cdot \times GL_{l_{r}}$ \\ 
$6:G=SO_{2n+1},Sp_{2n}$ or $SO_{2n}$ and $l_{r+1}=1$ & $L_{G,I}\simeq
GL_{l_{1}}\times \cdot \cdot \cdot \times GL_{l_{r}}\times SL_{2}$ \\ 
$7:G=SO_{2n}$ and $l_{r+1}=2$ & $L_{G,I}\simeq GL_{l_{1}}\times \cdot \cdot
\cdot \times GL_{l_{r}}\times SL_{2}\times SL_{2}$ \\ 
$8:G=SO_{2n}$ and $l_{r+1}=3$ & $L_{G,I}\simeq GL_{l_{1}}\times \cdot \cdot
\cdot \times GL_{l_{r}}\times SL_{4}$ \\ \hline
\end{tabular}
\label{dec}
\end{equation}
Note that $l_{r+1}=p$ means that $\alpha _{j_{r}}=\alpha _{n-p}.\;$The
factors of the decomposition of $L_{G,I}$ in a direct product of classical
groups are giving by the connected components of the diagram obtained by
removing the nodes corresponding to the simple roots $\alpha
_{j_{1}},...,\alpha _{j_{r}}$ in the Dynkin diagram of the root system of $%
\frak{g}.$

\noindent Since the Levi group $L_{G,I}$ is isomorphic to a direct product
of classical groups and contains the maximal torus $T_{G}$, we can define
the partition function $\mathcal{P}_{I}$ associated to $G_{0}=L_{G,I}$ as in
(\ref{deG0}) by the formal identity 
\begin{equation*}
\prod_{\alpha \in S_{G,I}}\frac{1}{1-e^{\alpha }}=\sum_{\beta \in \mathbb{Z}%
^{n}}\mathcal{P}^{G,I}(\beta )e^{\beta }
\end{equation*}
where $S_{G,I}=R_{G}^{+}-R_{G,I}^{+}$.\ Note that $S_{G,I}$ does not
coincide in general with the subset of positive roots of $R_{G}^{+}$
obtained as $\mathbb{N}$-linear combinations of the simple roots $\alpha
_{j_{1}},...,\alpha _{j_{r}}.\;$We describe in the table below, the sets $%
S_{G,I}$ corresponding to the decompositions of $L_{G,I}$ given in (\ref{dec}%
). Set 
\begin{eqnarray*}
\Theta _{G} &=&\left\{ 
\begin{tabular}{l}
$\{\varepsilon _{i}+\varepsilon _{j}\mid 1\leq i<j\leq n\}\cup \{\varepsilon
_{i}\mid 1\leq i\leq n\}$ if $G=SO_{2n+1}$ \\ 
$\{\varepsilon _{i}+\varepsilon _{j}\mid 1\leq i\leq j\leq n\}$ if $%
G=Sp_{2n+1}$ \\ 
$\{\varepsilon _{i}+\varepsilon _{j}\mid 1\leq i<j\leq n\}$ if $G=SO_{2n}$%
\end{tabular}
\right. \text{ and} \\
\Theta _{G}^{\ast } &=&\left\{ 
\begin{tabular}{l}
$\Theta _{G}-\{\varepsilon _{n}\}$ if $G=SO_{2n+1}$ \\ 
$\Theta _{G}-\{2\varepsilon _{n}\}$ if $G=Sp_{2n+1}$ \\ 
$\Theta _{G}-\{\varepsilon _{n-1}+\varepsilon _{n}\}$ if $G=SO_{2n}$%
\end{tabular}
\right. .
\end{eqnarray*}
\begin{equation*}
\begin{tabular}{|l|}
\hline
$1:S_{G,I}=\overset{r}{\underset{s=1}{\bigcup }}\{\varepsilon
_{i}-\varepsilon _{j}\mid 1\leq i\leq j_{s}<j\leq n\}$ \\ 
$2:S_{G,I}=\overset{r}{\underset{s=1}{\bigcup }}\{\varepsilon
_{i}-\varepsilon _{j}\mid 1\leq i\leq j_{s}<j\leq n\}\cup \{\varepsilon
_{i}+\varepsilon _{j}\mid 1\leq i<j\leq n$ and $i\leq j_{r}\}\cup
\{\varepsilon _{i}\mid 1\leq i\leq j_{r}\}$ \\ 
$3:S_{G,I}=\overset{r}{\underset{s=1}{\bigcup }}\{\varepsilon
_{i}-\varepsilon _{j}\mid 1\leq i\leq j_{s}<j\leq n\}\cup \{\varepsilon
_{i}+\varepsilon _{j}\mid 1\leq i\leq j\leq n$ and $i\leq j_{r}\}$ \\ 
$4:S_{G,I}=\overset{r}{\underset{s=1}{\bigcup }}\{\varepsilon
_{i}-\varepsilon _{j}\mid 1\leq i\leq j_{s}<j\leq n\}\cup \{\varepsilon
_{i}+\varepsilon _{j}\mid 1\leq i<j\leq n$ and $i\leq j_{r}\}$ \\ 
$5:S_{G,I}=\overset{r}{\underset{s=1}{\bigcup }}\{\varepsilon
_{i}-\varepsilon _{j}\mid 1\leq i\leq j_{s}<j\leq n\}\cup \Theta _{G}$ \\ 
$6:S_{G,I}=\overset{r}{\underset{s=1}{\bigcup }}\{\varepsilon
_{i}-\varepsilon _{j}\mid 1\leq i\leq j_{s}<j\leq n\}\cup \Theta _{G}^{\ast
} $ \\ 
$7:S_{G,I}=\overset{r}{\underset{s=1}{\bigcup }}\{\varepsilon
_{i}-\varepsilon _{j}\mid 1\leq i\leq j_{s}<j\leq n\}\cup \Theta _{G}^{\ast
} $ \\ 
$8:S_{G,I}=\overset{r}{\underset{s=1}{\bigcup }}\{\varepsilon
_{i}-\varepsilon _{j}\mid 1\leq i\leq j_{s}<j\leq n\}\cup \{\varepsilon
_{i}+\varepsilon _{j}\mid 1\leq i<j\leq n$ and $i\leq n-3\}$ \\ \hline
\end{tabular}
\end{equation*}

\newpage

\begin{example}
\label{exam1}Consider $G=Sp_{8}.\;$We give below the $16$ possible sets $I$
and $S_{G,I}$ for each Levi subgroup $L_{G,I}:$%
\begin{equation*}
\begin{tabular}{|ccc|}
\hline
$I$ & \multicolumn{1}{|c}{$L_{G,I}$} & \multicolumn{1}{|c|}{$S_{G,I}$} \\ 
\hline
$\{\alpha _{1},\alpha _{2},\alpha _{3},\alpha _{4}\}$ & \multicolumn{1}{|c|}{%
$Sp_{8}$} & $\emptyset $ \\ 
$\{\alpha _{2},\alpha _{3},\alpha _{4}\}$ & \multicolumn{1}{|c|}{$%
GL_{1}\times Sp_{6}$} & $\varepsilon _{1}\pm \varepsilon _{2},\varepsilon
_{1}\pm \varepsilon _{3},\varepsilon _{1}\pm \varepsilon _{3},2\varepsilon
_{1}$ \\ 
\multicolumn{1}{|c|}{$\{\alpha _{1},\alpha _{3},\alpha _{4}\}$} & 
\multicolumn{1}{c|}{$GL_{2}\times Sp_{4}$} & $
\begin{tabular}{c}
$\varepsilon _{1}\pm \varepsilon _{3},\varepsilon _{1}\pm \varepsilon
_{4},\varepsilon _{2}\pm \varepsilon _{3},\varepsilon _{2}\pm \varepsilon
_{4}$ \\ 
$\varepsilon _{1}+\varepsilon _{2},2\varepsilon _{1},2\varepsilon _{2}$%
\end{tabular}
$ \\ 
\multicolumn{1}{|c|}{$\{\alpha _{3},\alpha _{4}\}$} & \multicolumn{1}{c|}{$%
GL_{1}\times GL_{1}\times Sp_{4}$} & $
\begin{tabular}{c}
$\varepsilon _{1}\pm \varepsilon _{2},\varepsilon _{1}\pm \varepsilon
_{3},\varepsilon _{1}\pm \varepsilon _{4},\varepsilon _{2}\pm \varepsilon
_{3},\varepsilon _{2}\pm \varepsilon _{4}$ \\ 
$2\varepsilon _{1},2\varepsilon _{2}$%
\end{tabular}
$ \\ 
\multicolumn{1}{|c|}{$\{\alpha _{1},\alpha _{2},\alpha _{4}\}$} & 
\multicolumn{1}{c|}{$GL_{3}\times SL_{2}$} & $
\begin{tabular}{c}
$\varepsilon _{1}+\varepsilon _{2},\varepsilon _{1}+\varepsilon
_{3},\varepsilon _{1}\pm \varepsilon _{4},\varepsilon _{2}+\varepsilon
_{3},\varepsilon _{2}\pm \varepsilon _{4}$ \\ 
$\varepsilon _{3}\pm \varepsilon _{4},2\varepsilon _{1},2\varepsilon
_{2},2\varepsilon _{3}$%
\end{tabular}
$ \\ 
\multicolumn{1}{|c|}{$\{\alpha _{1},\alpha _{4}\}$} & \multicolumn{1}{c|}{$%
GL_{2}\times GL_{1}\times SL_{2}$} & $
\begin{tabular}{c}
$\varepsilon _{1}+\varepsilon _{2},\varepsilon _{1}\pm \varepsilon
_{3},\varepsilon _{1}\pm \varepsilon _{4},\varepsilon _{2}\pm \varepsilon
_{3},\varepsilon _{2}\pm \varepsilon _{4}$ \\ 
$\varepsilon _{3}\pm \varepsilon _{4},2\varepsilon _{1},2\varepsilon
_{2},2\varepsilon _{3}$%
\end{tabular}
$ \\ 
\multicolumn{1}{|c|}{$\{\alpha _{2},\alpha _{4}\}$} & \multicolumn{1}{c|}{$%
GL_{1}\times GL_{2}\times SL_{2}$} & $
\begin{tabular}{c}
$\varepsilon _{1}\pm \varepsilon _{2},\varepsilon _{1}\pm \varepsilon
_{3},\varepsilon _{1}\pm \varepsilon _{4},\varepsilon _{2}+\varepsilon
_{3},\varepsilon _{2}\pm \varepsilon _{4}$ \\ 
$\varepsilon _{3}\pm \varepsilon _{4},2\varepsilon _{1},2\varepsilon
_{2},2\varepsilon _{3}$%
\end{tabular}
$ \\ 
\multicolumn{1}{|c|}{$\{\alpha _{4}\}$} & \multicolumn{1}{c|}{$GL_{1}\times
GL_{1}\times GL_{1}\times SL_{2}$} & $
\begin{tabular}{c}
$\varepsilon _{1}\pm \varepsilon _{2},\varepsilon _{1}\pm \varepsilon
_{3},\varepsilon _{1}\pm \varepsilon _{4},\varepsilon _{2}\pm \varepsilon
_{3},\varepsilon _{2}\pm \varepsilon _{4}$ \\ 
$\varepsilon _{3}\pm \varepsilon _{4},2\varepsilon _{1},2\varepsilon
_{2},2\varepsilon _{3}$%
\end{tabular}
$ \\ 
\multicolumn{1}{|c|}{$\{\alpha _{1},\alpha _{2},\alpha _{3}\}$} & 
\multicolumn{1}{c|}{$GL_{4}$} & $
\begin{tabular}{c}
$\varepsilon _{1}+\varepsilon _{2},\varepsilon _{1}+\varepsilon
_{3},\varepsilon _{1}+\varepsilon _{4},\varepsilon _{2}+\varepsilon
_{3},\varepsilon _{2}+\varepsilon _{4}$ \\ 
$\varepsilon _{3}+\varepsilon _{4},2\varepsilon _{1},2\varepsilon
_{2},2\varepsilon _{3},2\varepsilon _{4}$%
\end{tabular}
$ \\ 
\multicolumn{1}{|c|}{$\{\alpha _{2},\alpha _{3}\}$} & \multicolumn{1}{c|}{$%
GL_{1}\times GL_{3}$} & $
\begin{tabular}{c}
$\varepsilon _{1}\pm \varepsilon _{2},\varepsilon _{1}\pm \varepsilon
_{3},\varepsilon _{1}\pm \varepsilon _{4},\varepsilon _{2}+\varepsilon
_{3},\varepsilon _{2}+\varepsilon _{4}$ \\ 
$\varepsilon _{3}+\varepsilon _{4},2\varepsilon _{1},2\varepsilon
_{2},2\varepsilon _{3},2\varepsilon _{4}$%
\end{tabular}
$ \\ 
\multicolumn{1}{|c|}{$\{\alpha _{1},\alpha _{2}\}$} & \multicolumn{1}{c|}{$%
GL_{3}\times GL_{1}$} & $
\begin{tabular}{c}
$\varepsilon _{1}+\varepsilon _{2},\varepsilon _{1}+\varepsilon
_{3},\varepsilon _{1}\pm \varepsilon _{4},\varepsilon _{2}+\varepsilon
_{3},\varepsilon _{2}\pm \varepsilon _{4}$ \\ 
$\varepsilon _{3}\pm \varepsilon _{4},2\varepsilon _{1},2\varepsilon
_{2},2\varepsilon _{3},2\varepsilon _{4}$%
\end{tabular}
$ \\ 
\multicolumn{1}{|c|}{$\{\alpha _{1},\alpha _{3}\}$} & \multicolumn{1}{c|}{$%
GL_{2}\times GL_{2}$} & $
\begin{tabular}{c}
$\varepsilon _{1}+\varepsilon _{2},\varepsilon _{1}\pm \varepsilon
_{3},\varepsilon _{1}\pm \varepsilon _{4},\varepsilon _{2}\pm \varepsilon
_{3},\varepsilon _{2}\pm \varepsilon _{4}$ \\ 
$\varepsilon _{3}+\varepsilon _{4},2\varepsilon _{1},2\varepsilon
_{2},2\varepsilon _{3},2\varepsilon _{4}$%
\end{tabular}
$ \\ 
\multicolumn{1}{|c|}{$\{\alpha _{3}\}$} & \multicolumn{1}{c|}{$GL_{1}\times
GL_{1}\times GL_{2}$} & $
\begin{tabular}{c}
$\varepsilon _{1}\pm \varepsilon _{2},\varepsilon _{1}\pm \varepsilon
_{3},\varepsilon _{1}\pm \varepsilon _{4},\varepsilon _{2}\pm \varepsilon
_{3},\varepsilon _{2}\pm \varepsilon _{4}$ \\ 
$\varepsilon _{3}+\varepsilon _{4},2\varepsilon _{1},2\varepsilon
_{2},2\varepsilon _{3},2\varepsilon _{4}$%
\end{tabular}
$ \\ 
\multicolumn{1}{|c|}{$\{\alpha _{2}\}$} & \multicolumn{1}{c|}{$GL_{1}\times
GL_{2}\times GL_{1}$} & $
\begin{tabular}{c}
$\varepsilon _{1}\pm \varepsilon _{2},\varepsilon _{1}\pm \varepsilon
_{3},\varepsilon _{1}\pm \varepsilon _{4},\varepsilon _{2}+\varepsilon
_{3},\varepsilon _{2}\pm \varepsilon _{4}$ \\ 
$\varepsilon _{3}\pm \varepsilon _{4},2\varepsilon _{1},2\varepsilon
_{2},2\varepsilon _{3},2\varepsilon _{4}$%
\end{tabular}
$ \\ 
\multicolumn{1}{|c|}{$\{\alpha _{1}\}$} & \multicolumn{1}{c|}{$GL_{2}\times
GL_{1}\times GL_{1}$} & $
\begin{tabular}{c}
$\varepsilon _{1}+\varepsilon _{2},\varepsilon _{1}\pm \varepsilon
_{3},\varepsilon _{1}\pm \varepsilon _{4},\varepsilon _{2}\pm \varepsilon
_{3},\varepsilon _{2}\pm \varepsilon _{4}$ \\ 
$\varepsilon _{3}\pm \varepsilon _{4},2\varepsilon _{1},2\varepsilon
_{2},2\varepsilon _{3},2\varepsilon _{4}$%
\end{tabular}
$ \\ 
\multicolumn{1}{|c|}{$\emptyset $} & \multicolumn{1}{c|}{$GL_{1}\times
GL_{1}\times GL_{1}\times GL_{1}$} & $
\begin{tabular}{c}
$\varepsilon _{1}\pm \varepsilon _{2},\varepsilon _{1}\pm \varepsilon
_{3},\varepsilon _{1}\pm \varepsilon _{4},\varepsilon _{2}\pm \varepsilon
_{3},\varepsilon _{2}\pm \varepsilon _{4}$ \\ 
$\varepsilon _{3}\pm \varepsilon _{4},2\varepsilon _{1},2\varepsilon
_{2},2\varepsilon _{3},2\varepsilon _{4}$%
\end{tabular}
$ \\ \hline
\end{tabular}
.
\end{equation*}
Note that $GL_{1}\times GL_{2}\times SL_{2}\simeq GL_{1}\times GL_{3}$ and $%
GL_{1}\times GL_{1}\times GL_{1}\times SL_{2}\simeq GL_{1}\times
GL_{1}\times GL_{2}.$
\end{example}

\noindent Suppose that the Levi subgroup $L_{G,I}$ decomposes following (\ref
{dec}) in 
\begin{equation}
L_{G,I}\simeq G_{1}\times \cdot \cdot \cdot \times G_{p}  \label{dec_L}
\end{equation}
where the $G_{k}$'s are classical Lie groups.\ Write $P_{G,I}^{+}$ for the
cone of dominant weights of $L_{G,I}.\;$The dominant weights of $P_{G,I}^{+}$
can be regarded as sequences $\mathbf{\mu }=(\mu ^{(1)},...,\mu ^{(p)})$
such that $\mu ^{(s)}$ belongs to $\widetilde{\mathcal{P}}_{l}$ if $%
G_{s}=GL_{l}$ and $\mu ^{(p)}$ belongs to $\mathcal{P}_{l}$ if $%
G_{p}=SL_{l},SO_{2l+1},Sp_{2l}$ or $SO_{2l}.$ Recall that, for such a
dominant weight $\mathbf{\mu }=(\mu ^{(1)},...,\mu ^{(p)})$, we denote by $%
\mu \in \mathbb{Z}^{n}$ the sequence obtained by reading successively the
parts of $\mu ^{(1)},...,\mu ^{(p)}$ from left to right.

\noindent We deduce immediately from Theorem \ref{th_branch} the branching
coefficients for the restriction of $V^{G}(\lambda )$ to the Levi subgroup $%
L_{G,I}.$

\begin{theorem}
\label{Th_B_levi}Consider $\lambda \in \mathcal{P}_{n}$ and $\mathbf{\mu }%
\in P_{G,I}^{+}$ then 
\begin{equation*}
\lbrack V(\lambda )^{G}:V(\mathbf{\mu })^{L_{G,I}}]=\sum_{w\in
W_{G}}(-1)^{\ell (w)}\mathcal{P}^{G,I}(w\circ \lambda -\mu ).
\end{equation*}
\end{theorem}

\noindent Note that $L_{G,I}\simeq GL_{n}$ when $I=\{\alpha _{1},...,\alpha
_{n-1}\}.\;$In this case, the branching coefficients $[V(\lambda
)^{G}:V(\gamma )^{GL_{n}}]$ where $\gamma \in \widetilde{\mathcal{P}}_{n}$
can be expressed in terms of the Littlewood-Richardson coefficients $%
c_{\gamma ,\lambda }^{\nu }$.\ For each classical group $SO_{2n+1},Sp_{2n}$
and $SO_{2n},$ set 
\begin{equation}
\prod_{\alpha \in \Theta _{G}}(1-e^{\alpha })^{-1}=\sum_{\beta \in \mathbb{N}%
^{n}}r_{G}(\beta )x^{\beta }.  \label{RG}
\end{equation}
Denote by $\mathcal{P}_{n}^{(2)}$ and $\mathcal{P}_{n}^{(1,1)}$ the sub-sets
of $\mathcal{P}_{n}$ containing respectively the partitions with even rows
and the partitions with even columns.

\begin{proposition}
\label{prop_mul_sum}Consider $\nu \in \mathcal{P}_{n}$ and $\lambda
=(\lambda ^{+},\lambda ^{-})\in \widetilde{\mathcal{P}}_{l}.$ Then

\begin{enumerate}
\item  $[V(\nu )^{SO_{2n+1}}:V(\lambda )^{GL_{n}}]=\sum_{w\in
W_{B_{n}}}(-1)^{\ell (w)}r_{SO_{2n+1}}(w\circ \nu -(\lambda ^{+},\lambda
^{-}))=\sum_{\gamma ,\delta \in \mathcal{P}_{n}}c_{\gamma ,\delta }^{\nu
}c_{\lambda ^{+},\lambda ^{-}}^{\delta },$

\item  $[V(\nu )^{Sp_{2n}}:V(\lambda )^{GL_{n}}]=\sum_{w\in
W_{C_{n}}}(-1)^{\ell (w)}r_{Sp_{2n}}(w\circ \nu -(\lambda ^{+},\lambda
^{-}))=\sum_{\gamma ,\delta \in \mathcal{P}_{n}^{(2)}}c_{\gamma ,\delta
}^{\nu }c_{\lambda ^{+},\lambda ^{-}}^{\delta },$

\item  $[V(\nu )^{SO_{2n}}:V(\lambda )^{GL_{n}}]=\sum_{w\in
W_{D_{n}}}(-1)^{\ell (w)}r_{SO_{2n}}(w\circ \nu -(\lambda ^{+},\lambda
^{-}))=\sum_{\gamma ,\delta \in \mathcal{P}_{n}^{(1,1)}}c_{\gamma ,\delta
}^{\nu }c_{\lambda ^{+},\lambda ^{-}}^{\delta }.$
\end{enumerate}
\end{proposition}

\begin{proof}
The right equalities of the Proposition are obtained by Theorem \ref
{Th_B_levi}. The left follow from a classical result by Littlewood (see \cite
{Li} appendix p\ 295).
\end{proof}

\bigskip

\noindent \textbf{Remark: }When $\lambda $ is a partition (that is $\lambda
^{+}=\lambda $ and $\lambda ^{-}=\emptyset $) we have by the above
proposition 
\begin{equation}
\left\{ 
\begin{tabular}{l}
$\lbrack V(\nu )^{SO_{2n+1}}:V(\lambda )^{GL_{n}}]=\sum_{\gamma \in \mathcal{%
P}_{n}}c_{\gamma ,\lambda }^{\nu }$ \\ 
$\lbrack V(\nu )^{SP_{2n}}:V(\lambda )^{GL_{n}}]=\sum_{\gamma \in \mathcal{P}%
_{n}^{(2)}}c_{\gamma ,\lambda }^{\nu }$ \\ 
$\lbrack V(\nu )^{SO_{2n}}:V(\lambda )^{GL_{n}}]=\sum_{\gamma \in \mathcal{P}%
_{n}^{(1,1)}}c_{\gamma ,\lambda }^{\nu }$%
\end{tabular}
\right. .  \label{express}
\end{equation}
In particular for $\kappa =(1,...,1)\in \mathbb{N}^{n}$ and any nonnegative
integer $k$ we obtain 
\begin{equation*}
\lbrack V(\nu +k\kappa )^{G}:V(\lambda +k\kappa )^{GL_{n}}]=[V(\nu
)^{G}:V(\lambda )^{GL_{n}}].
\end{equation*}

\subsection{Branching coefficients associated to an orthogonal decomposition
of the root system $D_{n}\label{subsec_lit}$}

Consider a $r$-tuple $\eta =(\eta _{1},...,\eta _{r})$ of positive integers
summing $n.$ We associate to $\eta $ the orthogonal decomposition $D_{\eta
}=D_{\eta _{1}}\cup \cdot \cdot \cdot \cup D_{\eta _{r}}$ of the root system 
$D_{n}$ such that for any $k=1,...,r$%
\begin{equation*}
D_{\eta _{k}}=\{\pm \varepsilon _{i}\pm \varepsilon _{j}\mid \eta
_{k-1}+1\leq i<j\leq \eta _{k}\}
\end{equation*}
with $\eta _{0}=1.$ Then $SO_{2n}$ contains a subgroup $SO_{\eta }$ such
that 
\begin{equation*}
SO_{\eta }\simeq SO_{2\eta _{1}}\times \cdot \cdot \cdot \times SO_{2\eta
_{r}.}
\end{equation*}
Note that $SO_{\eta }$ is not a Levi subgroup of $SO_{2n}.$ The dominant
weights of $SO_{\eta }$ are the $r$-tuple of partitions $\mathbf{\mu }=(\mu
^{(1)},...,\mu ^{(r)})$ such that $\mu ^{(k)}$ belongs to $\mathcal{P}_{\eta
_{k}}$ for any $k=1,...,r.$ Since $SO_{\eta }$ contains the maximal torus of 
$SO_{2n}$ we can apply Theorem \ref{th_branch} with $G_{0}=SO_{\eta }.$ The
corresponding partition function is defined by the formal identity 
\begin{equation}
\prod_{(i,j)\in E_{\eta }}(1-\frac{x_{i}}{x_{j}})^{-1}\prod_{(r,s)\in
E_{\eta }}(1-x_{r}x_{s})^{-1}=\sum_{\beta \in \mathbb{Z}^{n}}\mathcal{P}%
^{\eta }(\beta )e^{\beta }.  \label{par_So}
\end{equation}
where $E_{\eta }=\cup _{2\leq p\leq r}\{(i,j)\mid 1\leq i\leq \eta
_{1}+\cdot \cdot \cdot +\eta _{p-1}<j\leq n\}.$

\begin{proposition}
\label{prop_rest_diret_So}Consider $\lambda $ a partition and $\mathbf{\mu }$
a dominant weight of $SO_{\eta }.$ Then 
\begin{equation*}
\lbrack V(\lambda )^{SO_{2n}}:V(\mathbf{\mu })^{SO_{\eta }}]=\sum_{w\in
W_{D_{n}}}(-1)^{\ell (w)}\mathcal{P}^{\eta }(w\circ \lambda -\mu ).
\end{equation*}
\end{proposition}

\noindent \textbf{Remark: }Although it is possible to obtain similar
branching coefficients starting from orthogonal decompositions of the root
systems $B_{n}$ and $C_{n},$ we do not use them in the sequel.

\section{Generalization of Lusztig $q$-analogues}

\subsection{Quantization of the partition functions associated to \ a Levi
subgroup\label{subsec_quant}}

Consider a classical group $G$ and $I$ a subset of $\Sigma _{G}^{+}$.\ We
associated to the Levi subgroup $L_{G,I}$ the $q$-partition function $%
\mathcal{P}_{q}^{G,I}$ defined from the formal identity 
\begin{equation}
\prod_{\alpha \in S_{G,I}}\frac{1}{1-qe^{\alpha }}=\sum_{\beta \in \mathbb{Z}%
^{n}}\mathcal{P}_{q}^{G,I}(\beta )e^{\beta }.  \label{qpartLevi}
\end{equation}

\begin{definition}
\label{def_K}Let $\lambda $ be partition of $\mathcal{P}_{n}$ and $\mathbf{%
\mu }$ a weight of $L_{G,I}.$ We denote by $K_{\lambda ,\mu }^{G,I}(q)$ the
polynomial 
\begin{equation}
K_{\lambda ,\mu }^{G,I}(q)=\sum_{w\in W_{G}}(-1)^{\ell (w)}\mathcal{P}%
_{q}^{G,I}(w\circ \lambda -\mu ).  \label{def_Kdis}
\end{equation}
\end{definition}

\noindent \textbf{Remark: }Since $\mathcal{P}_{q}^{G,I}(\beta )=0$ for any $%
\beta \in \mathbb{Z}^{n}$ with $\left| \beta \right| <0,$ we have $%
K_{\lambda ,\mu }^{G,I}(q)\neq 0$ only if $\left| \lambda \right| \geq
\left| \mu \right| .$

\bigskip

\noindent When $\mathbf{\mu }$ is a dominant weight of $L_{G,I},$ we deduce
from Theorem \ref{Th_B_levi} that the polynomial $K_{\lambda ,\mu }^{G,I}(q)$
is a $q$-analogue of the branching coefficient $[V(\lambda )^{G}:V(\mathbf{%
\mu })^{L_{G,I}}].\;$When $I=\emptyset ,$ that is when $S_{G,I}$ contains
all the simple roots of $G,$ $L_{G,I}$ coincide with the maximal torus of $%
G, $ thus $P_{G,I}^{+}=P_{G}^{+}$ and $K_{\lambda ,\mu }^{G,I}(q)$ is the
Lusztig $q$-analogue associated to the weight $\mu $ in $V(\lambda )^{G}.\;$%
If we suppose that $\mu $ is a partition, it is known \cite{bry} that $%
K_{\lambda ,\mu }^{G,I}(q)$ has nonnegative integer coefficients.

\noindent The polynomials $K_{\lambda ,\mu }^{G,I}(q)$ can also be defined
from the Hilbert series of the Euler characteristic associated to certain
graded virtual $G$-modules $\chi _{\mathbf{\mu }}$.\ When $\mathbf{\mu }$ is
a dominant weight stable under the action of the Weyl group of $L_{G,I},$
Broer has proved in Theorem 2.2 of \cite{broer} that the higher cohomology
vanishes in the Euler characteristic associated to $\chi _{\mathbf{\mu }}$.\
In this case, by the use of the Borel-Weil-Bott Theorem, its graded formal
character has a nonnegative expansion on $\{e^{\lambda }\mid \lambda \in 
\mathcal{P}_{n}\}.\;$This implies in particular that $K_{\lambda ,\mu
}^{G,I}(q)$ has nonnegative coefficients. Note that the results of \cite
{broer} does not require that $G$ is a classical Lie group.\ In the context
of this article, the dominant weight $\mathbf{\mu }=(\mu ^{(1)},...,\mu
^{(r)})$ is stable under the action of the Weyl group of $L_{G,I}$ if and
only if the $\mu ^{(k)}$'s are rectangular partitions of decreasing heights
and $\mu ^{(r)}=0$ when $L_{G,I}$ is not a direct product of linear groups.
This yields to the following theorem:

\begin{theorem}
\label{TH_broer}(from \cite{broer}) Consider $L_{G,I}$ a Levi subgroup of
the classical lie group $G.$ Let $\lambda $ be a partition of $\mathcal{P}%
_{n}$ and $\mathbf{\mu }=(\mu ^{(1)},...,\mu ^{(r)})$ a dominant weight of $%
L_{G,I}$ such that the $\mu ^{(k)}$'s are rectangular partitions of
decreasing heights with $\mu ^{(r)}=0$ when $L_{G,I}$ is not a direct
product of linear groups.\ Then $K_{\lambda ,\mu }^{G,I}(q)$ has nonnegative
coefficients.
\end{theorem}

\noindent When $G=GL_{n},$ the polynomial $K_{\lambda ,\mu }^{GL_{n},I}(q)$
can also be interpreted as a $q$-analogue of the generalized
Littlewood-Richardson coefficient $c_{\mu ^{(1)},...,\mu ^{(r)}}^{\lambda }$
giving the multiplicity of $V(\lambda )^{GL_{n}}$ in $V^{GL_{n}}(\mu
^{(1)})\otimes \cdot \cdot \cdot \otimes V^{GL_{n}}(\mu ^{(r)}).\;$Suppose
that the $\mu ^{(k)}$'s are rectangular partitions and denote by $X_{\lambda
,\mathbf{\mu }}^{\emptyset }(q)$ the one-dimensional sum defined from the
affine $A_{n-1}^{(1)}$-crystal $B_{\mathbf{\mu }}$ associated to $\mathbf{%
\mu }$ and the partition $\lambda $ (\cite{HKOTY}). In fact $X_{\lambda ,%
\mathbf{\mu }}^{\emptyset }(q)$ is defined up to a power of $q$ depending on
the normalization of the energy function $H_{\mathbf{\mu }}$ chosen on the
vertices of $B_{\mathbf{\mu }}.$ By using a Morris-type recurrence formula
for the Poincar\'{e} polynomials and a combinatorial description of the
polynomials $K_{\lambda ,\mu }^{GL_{n},I}(q)$, Shimozono has obtained the
following theorem:

\begin{theorem}
\cite{Sh0}\label{th-shimo} Let $\lambda $ be a partition of $\mathcal{P}_{n}$
and $\mathbf{\mu }=(\mu ^{(1)},...,\mu ^{(r)})$ a dominant weight of $%
L_{G,I} $ such that the $\mu ^{(k)}$'s are rectangular partitions of
decreasing heights.\ Then 
\begin{equation*}
K_{\lambda ,\mu }^{GL_{n},I}(q)=q^{\ast }X_{\lambda ,\mathbf{\mu }%
}^{\emptyset }(q)
\end{equation*}
where $q^{\ast }$ is a power of $q$ depending on the normalization chosen
for $H_{\mathbf{\mu }}$.
\end{theorem}

\noindent \textbf{Remarks: }

\noindent $\mathrm{(i):}$ Theorem \ref{th-shimo} gives in particular a
combinatorial proof of the positivity of the polynomials $K_{\lambda ,\mu
}^{GL_{n},I}(q).\;$In the next paragraph we will use this result to derive
the positivity of the stable limits $\widetilde{K}_{\lambda ,\mu }^{G,I}(q)$
independently of Theorem \ref{TH_broer}.

\noindent $\mathrm{(ii):}$ Under the hypotheses of Theorem \ref{th-shimo},
it is conjectured that the polynomials $K_{\lambda ,\mu }^{GL_{n},I}(q)$
coincide with the $q$-analogues of the Littlewood-Richardson coefficients
introduced by Lascoux Leclerc and Thibon \cite{LLT}.

\bigskip

\noindent Numerous computations lead to conjecture that the positivity
result of Theorem \ref{TH_broer} can be extended to the case when the $r$%
-tuple $\mu $ associated to the dominant weight $\mathbf{\mu }\in
P_{G,I}^{+} $ is a partition.

\begin{conjecture}
\label{conj}Let $\lambda $ be partition of $\mathcal{P}_{n}$ and $\mathbf{%
\mu }$ a dominant weight of $L_{G,I}$ such that $\mu $ is a partition. Then $%
K_{\lambda ,\mu }^{G,I}(q)$ has nonnegative coefficients.
\end{conjecture}

\noindent I was informed that an equivalent statement of this conjecture
appeared for the first time in unpublished notes by Broer. In the
terminology of \cite{broer}, it is indeed equivalent to say that higher
cohomology vanishes in $\chi _{\mathbf{\mu }}$ when $\mu $ is a dominant
weight.

\begin{example}
Continuing Example \ref{exam1} with $G=Sp_{8},$ $\lambda =(4,2,2,1)$ and $%
\mu =(3,1,1,0),$ we obtain the following polynomials $K_{\lambda ,\mu
}^{G,I}(q)$%
\begin{equation*}
\begin{tabular}{|c|c|}
\hline
$L_{G,I}$ & $K_{\lambda ,\mu }^{G,I}(q)$%
\begin{tabular}{l}
\  \\ 
\ 
\end{tabular}
\\ \hline
$Sp_{8}$ & $0$ \\ 
$GL_{1}\times Sp_{6}$ & $0$ \\ 
$GL_{2}\times Sp_{4}$ & $2q^{4}$ \\ 
$GL_{1}\times GL_{1}\times Sp_{4}$ & $q^{3}+2q^{2}$ \\ 
$GL_{3}\times SL_{2}$ & $q^{3}+q^{2}$ \\ 
$GL_{2}\times GL_{1}\times SL_{2}$ & $3q^{4}+4q^{3}+q^{2}$ \\ 
$GL_{1}\times GL_{2}\times SL_{2}$ & $q^{4}+2q^{3}+q^{2}$ \\ 
$GL_{1}\times GL_{1}\times GL_{1}\times SL_{2}$ & $%
2q^{5}+4q^{4}+4q^{3}+q^{2} $ \\ 
$GL_{4}$ & $q^{2}$ \\ 
$GL_{1}\times GL_{3}$ & $q^{3}+q^{2}$ \\ 
$GL_{3}\times GL_{1}$ & $q^{3}+2q^{3}+q^{2}$ \\ 
$GL_{2}\times GL_{2}$ & $3q^{4}+4q^{3}+q^{2}$ \\ 
$GL_{1}\times GL_{1}\times GL_{2}$ & $2q^{5}+4q^{4}+4q^{3}+q^{2}$ \\ 
$GL_{1}\times GL_{2}\times GL_{1}$ & $q^{5}+2q^{4}+3q^{2}+q^{2}$ \\ 
$GL_{2}\times GL_{1}\times GL_{1}$ & $%
q^{7}+2q^{6}+3q^{5}+4q^{4}+4q^{3}+q^{2} $ \\ 
$GL_{1}\times GL_{1}\times GL_{1}\times GL_{1}$ & $%
q^{8}+2q^{7}+3q^{6}+4q^{5}+5q^{4}+4q^{3}+q^{2}$ \\ \hline
\end{tabular}
.
\end{equation*}
When $L_{G,I}=GL_{1}\times GL_{1}\times GL_{1}\times GL_{1},$ $K_{\lambda
,\mu }^{G,I}(q)$ is the Lusztig $q$-analogue for the root system $C_{4}$
associated to the partitions $\lambda $ and $\mu .$ Note also that the
polynomials corresponding to the isomorphic Levi subgroups $GL_{1}\times
GL_{2}\times SL_{2}\simeq GL_{1}\times GL_{3}$ and $GL_{1}\times
GL_{1}\times GL_{1}\times SL_{2}\simeq GL_{1}\times GL_{1}\times GL_{2}$ are
equal.
\end{example}

\bigskip

\noindent In the sequel we will also led to consider another family of $q$%
-analogues for the branching coefficients $[V(\lambda )^{SO_{2n+1}}:V(%
\mathbf{\mu })^{L,SO_{2n+1}}]$ obtained from the partition function $%
\mathcal{P}_{q,h}^{SO_{2n+1},I}$ defined by the expansion 
\begin{equation*}
\prod_{\alpha \in S_{SO_{2n+1},I}}\frac{1}{1-q^{h(\alpha )}e^{\alpha }}%
=\sum_{\beta \in \mathbb{Z}^{n}}\mathcal{P}_{q,h}^{SO_{2n+1},I}(\beta
)e^{\beta }
\end{equation*}
where $h(\alpha )=2$ if $\alpha =\varepsilon _{i},$ $i=1,...,n$ and $%
h(\alpha )=1$ otherwise.

\begin{definition}
Let $\lambda $ be partition of $\mathcal{P}_{n}$ and $\mathbf{\mu }$ a
dominant weight of $L_{SO_{2n+1},I}.$ We denote by $\mathcal{K}_{\lambda
,\mu }^{SO_{2n+1},I}(q)$ the polynomial 
\begin{equation*}
\mathcal{K}_{\lambda ,\mu }^{SO_{2n+1},I}(q)=\sum_{w\in W_{B_{n}}}(-1)^{\ell
(w)}\mathcal{P}_{q,h}^{SO_{2n+1},I}(w\circ \lambda -\mu ).
\end{equation*}
\end{definition}

\noindent We have clearly $\mathcal{K}_{\lambda ,\mu
}^{SO_{2n+1},I}(1)=K_{\lambda ,\mu }^{SO_{2n+1},I}(1).\;$Note that the
polynomials $\mathcal{K}_{\lambda ,\mu }^{SO_{2n+1},I}(q)$ can have negative
coefficients even if the hypotheses of Conjecture \ref{conj} are verified.\
Nevertheless we are going to see that they admit a stable limit which
decompose as nonnegative combination of Poincar\'{e} polynomials.\ This is
not the case for the polynomials $K_{\lambda ,\mu }^{SO_{2n+1},I}(q)$ which
implies that they can not coincide with one-dimensional sums (see Conjecture 
\ref{conj_last}).$\;$This situation is analogous to that observed in \cite
{LS} where the one-dimensional sums considered are, for affine crystals of
kind $(1),$ equal to Lusztig $q$-analogues related to affine Hecke algebras
of type $B_{n}$ with parameters $q$ and $q^{2}$.

\subsection{Stable limit}

With the notation of the above paragraph, we define the stable limit $%
\widetilde{K}_{\lambda ,\mu }^{G,I}(q)$ by setting 
\begin{equation}
\widetilde{K}_{\lambda ,\mu }^{G,I}(q)=\left\{ 
\begin{tabular}{l}
$\sum_{\sigma \in S_{n}}(-1)^{\ell (\sigma )}\mathcal{P}_{q}^{G,I}(\sigma
\circ \lambda -\mu )$ if $G=GL_{n},Sp_{2n}$ or $SO_{2n}$ \\ 
$\sum_{\sigma \in S_{n}}(-1)^{\ell (\sigma )}\mathcal{P}%
_{q,h}^{SO_{2n+1},I}(\sigma \circ \lambda -\mu )$ if $G=SO_{2n+1}$%
\end{tabular}
\right. .  \label{def_Ktilde}
\end{equation}
This is an expression for $G$ a classical Lie group, but the sum runs over
the parabolic subgroup of $W_{G}$ generated by the $s_{i},$ $i=1,...,n$
which is a copy of the Weyl group of the root system $A_{n-1}.$ We have $%
\widetilde{K}_{\lambda ,\mathbf{\mu }}^{GL_{n},I}(q)=K_{\lambda ,\mathbf{\mu 
}}^{GL_{n},I}(q)$ since the Weyl group of $GL_{n}$ is the symmetric group $%
S_{n}.\;$

\begin{lemma}
\label{lem_ktilde}Consider $\lambda ,\mu $ two partitions of length $n$ such
that $\left| \lambda \right| \geq \left| \mu \right| .\;$Let $k$ be any
integer such that $k\geq \frac{\left| \lambda \right| -\left| \mu \right| }{2%
}$.\ Then we have 
\begin{equation*}
\widetilde{K}_{\lambda ,\mu }^{G,I}(q)=\left\{ 
\begin{tabular}{l}
$K_{\lambda +k\kappa ,\mathbf{\mu }+k\kappa }^{G,I}(q)$ when $%
G=GL_{n},Sp_{2n}$ or $SO_{2n}$ \\ 
$\mathcal{K}_{\lambda +k\kappa ,\mathbf{\mu }+k\kappa }^{SO_{2n+1},I}(q)$
when $G=SO_{2n+1}$%
\end{tabular}
\right. .
\end{equation*}
where $\kappa =(1,...,1)\in \mathbb{N}^{n}.$
\end{lemma}

\begin{proof}
Suppose $G=GL_{n},Sp_{2n}$ or $SO_{2n}.\;$For any $\beta \in \mathbb{Z}^{n},$
we have $\mathcal{P}_{q}^{G,I}(\beta )=0$ if $\beta $ is not a linear
combination of the positive roots of $S_{G,I}$ with nonnegative
coefficients. This implies that $\mathcal{P}_{q}^{G,I}(\beta )=0$ if $\left|
\beta \right| <0.$ Then the Lemma is a consequence of Lemma \ref{lem_tech}
applied with $\mathcal{M}=\mathcal{P}_{q}^{G,I}$. We proceed similarly for $%
G=SO_{2n+1}$ by using $\mathcal{P}_{q}^{G,I}$ instead of $\mathcal{P}%
_{q,h}^{SO_{2n+1},I}$.
\end{proof}

\ \bigskip

\noindent \textbf{Remark: }Since $\sigma (\kappa )=\kappa $ for any $\sigma
\in S_{n},$ we have the following stability property 
\begin{equation*}
\widetilde{K}_{\lambda +k\kappa ,\mathbf{\mu }+k\kappa }^{G,I}(q)=\widetilde{%
K}_{\lambda ,\mu }^{G,I}(q)\text{ for any integer }k
\end{equation*}
which justifies the above terminology.\ So we can extend the definition of $%
\widetilde{K}_{\lambda ,\mu }^{G,I}(q)$ when $\lambda $ and $\mu $ are
decreasing sequences of integers (positive or not).

\bigskip

\noindent For any $\xi \in \mathbb{Z}^{n},$ we define the polynomial $K_{\xi
,\mu }^{GL_{n},I}(q)$ by replacing, in (\ref{def_Kdis}) the partition $%
\lambda $ by $\xi .\;$There exists a straightening procedure for the
polynomials $K_{\xi ,\mu }^{GL_{n},I}(q)$ which follows immediately from the
fact that the set $\{\sigma \circ \xi \mid \sigma \in S_{n}\}$ (see (\ref
{dotaction})) intersects at most one time the cone of dominant weights of $%
GL_{n}$.

\begin{lemma}
\label{strilaw}Consider $\mu $ and $\xi $ in $\mathbb{Z}^{n}.\;$Then 
\begin{equation*}
\left\{ 
\begin{tabular}{l}
$K_{\xi ,\mu }^{GL_{n},I}(q)=(-1)^{l(\tau )}K_{\nu ,\mu }^{GL_{n},I}(q)$ if $%
\xi =\tau \circ (\nu )$ with $\tau \in S_{n}$ and $\nu \in \widetilde{%
\mathcal{P}}_{n}$ \\ 
$0$ otherwise
\end{tabular}
\right.
\end{equation*}
where $\widetilde{\mathcal{P}}_{n}=\{\gamma =(\gamma _{1},...,\gamma
_{n})\in \mathbb{Z}^{n},\gamma _{1}\geq \gamma _{2}\geq \cdot \cdot \cdot
\geq \gamma _{n}\}$.
\end{lemma}

\bigskip

In the sequel of this paragraph we restrict ourselves to the case when $%
G=SO_{2n+1},Sp_{2n}$ or $SO_{2n}$ and $l_{r+1}=0$ (with the notation of
paragraph \ref{subsec_levi}).$\;$This corresponds to the decomposition $5$
given in table (\ref{dec}), that is we suppose that $I$ does not contains
the simple root $\alpha _{n}$. In this case 
\begin{equation*}
L_{G,I}\simeq GL_{l_{1}}\times \cdot \cdot \cdot \times GL_{l_{r}}.
\end{equation*}
and $I$ is also a subset of $\Sigma _{GL_{n}}^{+},$ thus determinates a Levi
subgroup of $GL_{n}$ which is isomorphic to $L_{G,I}$.\ Moreover we have $%
S_{G,I}=S_{GL_{n},I}\cup \Theta _{G}.$ The $q$-partition functions $\mathcal{%
P}_{q}^{G,I},$ $G=Sp_{2n}$ or $SO_{2n}$ and $\mathcal{P}_{q,h}^{SO_{2n+1},I}$
can be expressed in terms of the $q$-partition $\mathcal{P}_{q}^{GL_{n},I}$:

\begin{lemma}
\label{lem_util}For any $\beta \in \mathbb{Z}^{n}$ we have 
\begin{equation*}
\left\{ 
\begin{tabular}{l}
$\mathrm{(i):}$ $\mathcal{P}_{q}^{G,I}(\beta )=q^{\left| \beta \right|
/2}\sum_{\delta \in \mathbb{N}^{n},\left| \delta \right| =\left| \beta
\right| }r_{G}(\delta )\mathcal{P}_{q}^{GL_{n},I}(\beta -\delta )\text{ for }%
G=Sp_{2n}\text{ or }SO_{2n}$ \\ 
$\mathrm{(ii):}$ $\mathcal{P}_{q,h}^{SO_{2n+1},I}(\beta )=q^{\left| \beta
\right| /2}\sum_{\delta \in \mathbb{N}^{n},\left| \delta \right| =\left|
\beta \right| }r_{SO_{2n+1}}(\delta )\mathcal{P}_{q}^{GL_{n},I}(\beta
-\delta )$%
\end{tabular}
\right. .
\end{equation*}
\end{lemma}

\begin{proof}
$\mathrm{(i):}$ The $q$-partition function $\mathcal{P}_{q}^{GL_{n},I}$ is
defined by 
\begin{equation*}
\prod_{\alpha \in S_{GL_{n}},I}\left( 1-qe^{\alpha }\right)
^{-1}=\sum_{\gamma \in \mathbb{Z}^{n}}\mathcal{P}_{q}^{GL_{n},I}(\gamma
)e^{\gamma }
\end{equation*}
and since $S_{G,I}=S_{GL_{n},I}\cup \Theta _{G}$ the $q$-partition function $%
\mathcal{P}_{q}^{G,I}$ verifies 
\begin{equation*}
\sum_{\beta \in \mathbb{Z}^{n}}\mathcal{P}_{q}^{G,I}(\beta )e^{\beta
}=\prod_{\alpha \in \Theta _{G}}(1-qe^{\alpha })^{-1}\prod_{\alpha \in
S_{GL_{n}},I}\left( 1-qe^{\alpha }\right) ^{-1}
\end{equation*}
and we derive from (\ref{RG}) 
\begin{equation}
\prod_{\alpha \in \Theta _{G}}(1-qe^{\alpha })^{-1}=\sum_{\delta \in \mathbb{%
N}^{n}}q^{\left| \delta \right| /2}r_{G}(\delta )e^{\xi }  \label{equal_fund}
\end{equation}
since the number of roots appearing in a decomposition of $\delta \in 
\mathbb{N}^{n}$ as a sum of positive roots $\varepsilon _{r}+\varepsilon
_{s} $ with $1\leq r<s\leq n$ or $2\varepsilon _{i}$ with $1\leq i\leq n$ is
always equal to $\left| \delta \right| /2.$ Thus we obtain 
\begin{equation*}
\sum_{\beta \in \mathbb{Z}^{n}}\mathcal{P}_{q}^{G,I}(\beta )e^{\beta
}=\sum_{\gamma \in \mathbb{Z}^{n}}\sum_{\delta \in \mathbb{N}^{n}}q^{\left|
\delta \right| /2}\mathcal{P}_{q}^{GL_{n},I}(\gamma )r_{G}(\delta )e^{\delta
+\gamma }.
\end{equation*}
We derive the equality $\mathcal{P}_{q}^{G,I}(\beta )=\sum_{\gamma +\delta
=\beta }r_{G}(\delta )q^{\left| \delta \right| /2}\mathcal{P}%
_{q}^{GL_{n},I}(\gamma ).\;$Since the set $S_{GL_{n},I}$ contains only
positive roots $\alpha $ with $\left| \alpha \right| =0$, we will have $%
\mathcal{P}_{q}^{GL_{n},I}(\gamma )=0$ when $\left| \gamma \right| \neq 0.\;$%
So we can suppose $\left| \gamma \right| =0$ and $\left| \delta \right|
=\left| \beta \right| $ in the previous sum.

\noindent $\mathrm{(ii):}$ Since $h(\alpha )=2$ when $\left| \alpha \right|
=1$ we can also write 
\begin{equation*}
\prod_{\alpha \in \Theta _{SO_{2n+1}}}(1-q^{h(\alpha )}e^{\alpha
})^{-1}=\sum_{\delta \in \mathbb{N}^{n}}q^{\left| \delta \right|
/2}r_{SO_{2n+1}}(\delta )e^{\xi }.
\end{equation*}
Then we derive $\mathrm{(ii)}$ by proceeding as in $\mathrm{(i)}$.
\end{proof}

\noindent \textbf{Remark: }A similar result for the $q$-partition function $%
\mathcal{P}_{q}^{SO_{2n+1},I}$ does not exit. Indeed the number of roots
appearing in a decomposition of $\delta \in \mathbb{N}^{n}$ as a sum of
positive roots $\varepsilon _{r}+\varepsilon _{s}$ with $1\leq r<s\leq n$
and $\varepsilon _{i}$ with $1\leq i\leq n$ does not depend only of $\left|
\delta \right| $ since $\left| \varepsilon _{r}+\varepsilon _{s}\right| \neq
\left| \varepsilon _{i}\right| .$

\begin{theorem}
\label{prop_dec_K_c}Suppose $G=SO_{2n+1},Sp_{2n}$ or $SO_{2n}$ and $%
l_{r+1}=0 $ and consider $\lambda ,\mu \in \mathcal{P}_{n}$ such that $%
\left| \lambda \right| \geq \left| \mu \right| .$ Then for any integer $%
k\geq \frac{\left| \lambda \right| -\left| \mu \right| }{2}$we have: 
\begin{equation*}
\widetilde{K}_{\lambda ,\mu }^{G,I}(q)=q^{\tfrac{\left| \lambda \right|
-\left| \mu \right| }{2}}\sum_{\gamma \in \widetilde{\mathcal{P}}%
_{n}}[V(\lambda +k\kappa )^{G}:V(\gamma +k\kappa )^{GL_{n}}]K_{\gamma ,\mu
}^{GL_{n},I}(q).
\end{equation*}
\end{theorem}

\begin{proof}
We only give the proof for $G=Sp_{2n}$ or $SO_{2n}$, the case $G=SO_{2n+1}$
is similar.\ We have 
\begin{equation*}
\widetilde{K}_{\lambda ,\mu }^{G,I}(q)=\sum_{\sigma \in \mathcal{S}%
_{n}}(-1)^{\ell (\sigma )}\mathcal{P}_{q}^{G,I}(\sigma (\lambda +\rho )-(\mu
+\rho )).
\end{equation*}
Hence from the previous lemma we derive 
\begin{equation*}
\widetilde{K}_{\lambda ,\mu }^{G,I}(q)=\sum_{\sigma \in \mathcal{S}%
_{n}}(-1)^{\ell (\sigma )}\sum_{\delta \in \mathbb{N}^{n},\left| \delta
\right| =\left| \beta \right| }r_{G}(\delta )q^{\left| \beta \right| /2}%
\mathcal{P}_{q}^{GL_{n},I}(\sigma (\lambda +\rho )-(\mu +\delta +\rho ))
\end{equation*}
where $\beta =\sigma (\lambda +\rho )-(\mu +\rho )$ in the second sum.\
Since $\left| \beta \right| =\left| \lambda \right| -\left| \mu \right| ,$
we obtain 
\begin{equation*}
\widetilde{K}_{\lambda ,\mu }^{G,I}(q)=q^{\frac{\left| \lambda \right|
-\left| \mu \right| }{2}}\sum_{\sigma \in \mathcal{S}_{n}}(-1)^{\ell (\sigma
)}\sum_{\delta \in \mathbb{N}^{n},\left| \delta \right| =\left| \lambda
\right| -\left| \mu \right| }r_{G}(\delta )\mathcal{P}_{q}^{GL_{n},I}(\sigma
(\lambda +\rho -\sigma ^{-1}(\delta ))-(\mu +\rho ))
\end{equation*}
For any $\sigma \in \mathcal{S}_{n},$ we have $\sigma ^{-1}(\mathbb{N}^{n})=%
\mathbb{N}^{n}$ and $r_{G}(\delta )=r_{G}(\sigma (\delta ))$ since $\sigma
(\Theta _{G})=\Theta _{G}.$ Thus 
\begin{multline}
\widetilde{K}_{\lambda ,\mu }^{G,I}(q)=q^{\frac{\left| \lambda \right|
-\left| \mu \right| }{2}}\sum_{\sigma \in \mathcal{S}_{n}}(-1)^{\ell (\sigma
)}\sum_{\delta \in \mathbb{N}^{n},\left| \delta \right| =\left| \lambda
\right| -\left| \mu \right| }r_{G}(\delta )\mathcal{P}_{q}^{GL_{n},I}(\sigma
(\lambda +\rho -\delta )-(\mu +\rho ))=  \label{equK} \\
q^{\frac{\left| \lambda \right| -\left| \mu \right| }{2}}\sum_{\delta \in 
\mathbb{N}^{n},\left| \delta \right| =\left| \lambda \right| -\left| \mu
\right| }r_{G}(\delta )K_{\lambda -\delta ,\mu }^{GL_{n},I}(q).
\end{multline}
Now by Lemma \ref{strilaw}, $K_{\lambda -\delta ,\mu }^{GL_{n},I}(q)=0$ or
there exits $\sigma \in \mathcal{S}_{n}$ and $\gamma \in \widetilde{\mathcal{%
P}}_{n}$ such that $\gamma =\sigma ^{-1}\circ (\lambda -\delta ).$ Then we
have $\left| \gamma \right| =\left| \lambda \right| -\left| \delta \right|
=\left| \mu \right| $ and $\delta =\lambda +\rho -\sigma (\gamma +\rho )$.\
It follows that 
\begin{equation*}
\widetilde{K}_{\lambda ,\mu }^{G,I}(q)=q^{\frac{\left| \lambda \right|
-\left| \mu \right| }{2}}\sum_{\sigma \in \mathcal{S}_{n}}(-1)^{\ell (\sigma
)}\sum_{\gamma \in \widetilde{\mathcal{P}}_{n}}r_{G}(\lambda +\rho -\sigma
(\gamma +\rho ))K_{\gamma ,\mu }^{GL_{n},I}(q).
\end{equation*}
Since $c(\delta )=c(\sigma (\delta ))$ for any $\sigma \in \mathcal{S}_{n}$
and $\delta \in \mathbb{N}^{n}$, we obtain the equality 
\begin{equation*}
\widetilde{K}_{\lambda ,\mu }^{G,I}(q)=q^{\frac{\left| \lambda \right|
-\left| \mu \right| }{2}}\sum_{\gamma \in \widetilde{\mathcal{P}}%
_{n}}\sum_{\sigma \in \mathcal{S}_{n}}(-1)^{\ell (\sigma )}r_{G}(\sigma
(\lambda +\rho )-(\gamma +\rho ))K_{\gamma ,\mu }^{GL_{n},I}(q).
\end{equation*}
We deduce from Lemma \ref{lem_tech} applied with $\mathcal{M}=r_{G}$ and
from Proposition \ref{prop_mul_sum} that the equality 
\begin{equation*}
\sum_{\sigma \in \mathcal{S}_{n}}(-1)^{\ell (\sigma )}r_{G}(\sigma (\lambda
+\rho )-(\gamma +\rho ))=[V(\lambda +k\kappa )^{G}:V(\gamma +k\kappa
)^{GL_{n}}]
\end{equation*}
holds for any integer $k\geq \frac{\left| \lambda \right| -\left| \gamma
\right| }{2}=\frac{\left| \lambda \right| -\left| \mu \right| }{2}.$ This
yields to the desired equality.
\end{proof}

\bigskip

\noindent By using Theorem \ref{th-shimo} we obtain immediately

\begin{corollary}
\label{corKC}Suppose $G=SO_{2n+1},Sp_{2n}$ or $SO_{2n}$ and $l_{r+1}=0$.
Consider $\lambda \in \mathcal{P}_{n},$ $\mathbf{\mu }=(\mu ^{(1)},...,\mu
^{(p)})\in P_{G,I}^{+}$ such that the $\mu ^{(k)}$'s are rectangular
partitions of decreasing heights.\ Then $\widetilde{K}_{\lambda ,\mu
}^{G,I}(q)$ have nonnegative coefficients.
\end{corollary}

\noindent \textbf{Remark: }When\textbf{\ }$G=Sp_{2n}$ or $SO_{2n},$ Lemma 
\ref{lem_ktilde} implies that the above corollary can be regarded as a
particular case of Theorem \ref{TH_broer}. Since the polynomials $\widetilde{%
K}_{\lambda ,\mu }^{SO_{2n+1},I}(q)$ are generalized Lusztig $q$-analogues
defined by using the parameters $q$ and $q^{2},$ we can not deduce their
positivity from the results of Broer.

\section{Some dualities between tensor product and branching coefficients}

\subsection{Determinantal identities and operators on formal series\label%
{subsec_Jacobi-trudi}}

Consider $k,m\in \mathbb{Z}$ such that $m>0$. When $k$ is a nonnegative
integer, write $(k)_{n}=(k,0,...,0)$ for the partition of length $n$ with a
unique non-zero part equal to $k$.\ Then set $h_{k}^{G}=s_{(k)_{n}}^{G}$ if $%
k\geq 0$ and $h_{k}^{G}=0$ otherwise where $s_{(k)_{n}}^{G}$ is the
universal character of Koike and Terada associated to $(k)_{n}$ for the Lie
group $G$. Given $\alpha =(\alpha _{1},...,\alpha _{m})\in \mathbb{Z}^{m}$
define 
\begin{equation}
u_{\alpha }^{G}=\det \left( 
\begin{array}{cccc}
h_{\alpha _{1}}^{G} & h_{\alpha _{1}+1}^{G}+h_{\alpha _{1}-1}^{G} & \cdot
\cdot \cdot \cdot \cdot \cdot \cdot \cdot \cdot \cdot & h_{\alpha
_{1}+m-1}^{G}+h_{\alpha _{1}-m+1}^{G} \\ 
h_{\alpha _{2}-1}^{G} & h_{\alpha _{2}}^{G}+h_{\alpha _{2}-2}^{G} & \cdot
\cdot \cdot \cdot \cdot \cdot \cdot \cdot \cdot \cdot & h_{\alpha
_{2}+m-2}^{G}+h_{\alpha _{2}-m}^{G} \\ 
\cdot & \cdot & \cdot \cdot \cdot \cdot \cdot \cdot \cdot \cdot \cdot \cdot
& \cdot \\ 
\cdot & \cdot & \cdot \cdot \cdot \cdot \cdot \cdot \cdot \cdot \cdot \cdot
& \cdot \\ 
h_{\alpha _{m}-m+1}^{G} & h_{\alpha _{m}-m+2}^{G}+h_{\alpha _{n}-m}^{G} & 
\cdot \cdot \cdot \cdot \cdot \cdot \cdot \cdot \cdot \cdot & h_{\alpha
_{m}}^{G}+h_{\alpha _{m}-2m+2}^{G}
\end{array}
\right) .  \label{u_h}
\end{equation}
The following proposition is a well known analogue of the Jacobi-Trudi
determinantal formula for $G=SO_{2n+1},Sp_{2n}$ or $SO_{2n}.$

\begin{proposition}
\label{prop_fult}(see \cite{FH} \S 24.2) Consider $\lambda $ a partition
with at most $m$ nonzero parts.\ Then for $G=SO_{2n+1},Sp_{2n}$ or $SO_{2n}$
we have $u_{\lambda }^{G}=s_{\lambda }^{G}$.
\end{proposition}

\noindent By using elementary permutations on rows in the determinant (\ref
{u_h}) we obtain the straightening law for $u_{\alpha }^{G}$:

\begin{lemma}
\label{straight}Consider $\alpha \in \mathbb{Z}^{m}$ then 
\begin{equation*}
u_{\alpha }^{G}=\left\{ 
\begin{tabular}{l}
$(-1)^{\ell (\sigma )}s_{\lambda }^{G}$ if there exists $\sigma \in \mathcal{%
S}_{m}$ and $\lambda \in \mathcal{P}_{m}$ such that $\sigma \circ \alpha
=\lambda $ \\ 
$0$ otherwise
\end{tabular}
\right. .
\end{equation*}
\end{lemma}

\bigskip

\noindent Denote by $\mathcal{L}_{n}=\mathbb{K[}%
[x_{1},x_{1}^{-1},...,x_{n},x_{n}^{-1}]]$ the vector space of formal Laurent
series in the indeterminates $x_{1},x_{1}^{-1},...,x_{n},x_{n}^{-1}.\;$We
identify the ring of Laurent polynomials $L_{n}=\mathbb{K[}%
x_{1},x_{1}^{-1},...,x_{n},x_{n}^{-1}]$ with the sub-space of $\mathcal{L}%
_{n}$ containing the finite formal series.\ The vector space $\mathcal{L}%
_{n} $ is not a ring since the formal series are in the two directions. More
precisely, the product $F_{1}\cdot \cdot \cdot F_{r}$ of the formal series $%
F_{i}=\sum_{\beta _{i}\in E_{i}}x^{\beta _{i}}$ $i=1,...,r$ is defined if
and only if, for any $\gamma \in \mathbb{Z}^{n},$ the number $N_{\gamma }$
of decompositions $\gamma =\beta _{1}+\cdot \cdot \cdot +\beta _{r}$ such
that $\beta _{i\in E_{i}}$ is finite and in this case we have 
\begin{equation*}
F_{1}\cdot \cdot \cdot F_{r}=\sum_{\gamma \in \mathbb{Z}^{n}}N_{\gamma
}x^{\gamma }.
\end{equation*}
In particular the product $P\cdot F$ with $P\in L_{n}$ and $F\in \mathcal{L}%
_{n}$ is well defined.

\noindent Consider $\alpha \in \mathbb{Z}^{m}.$ We set $h_{\alpha
}^{G}=h_{\alpha _{1}}^{G}\cdot \cdot \cdot h_{\alpha _{m}}^{G}.$ Let $%
K=[k_{1},...,k_{m}]\subset \{1,...,n\}$ be the interval containing the $m$
consecutive integers $k_{1}<\cdot \cdot \cdot <k_{m}.$\ Denote by $\mathcal{L%
}_{K}\subset \mathcal{L}_{n}$ the vector space of formal Laurent formal
series in the indeterminates $%
x_{k_{1}},x_{k_{1}}^{-1},...,x_{k_{m}},x_{k_{m}}^{-1}.\;$We define the
determinant 
\begin{equation*}
\delta _{K}(\alpha )=\det \left( 
\begin{array}{cccc}
x_{k_{1}}^{\alpha _{1}} & x_{k_{1}}^{\alpha _{1}+1}+x_{k_{1}}^{\alpha _{1}-1}
& \cdot \cdot \cdot \cdot \cdot \cdot \cdot \cdot \cdot \cdot & 
x_{k_{1}}^{\alpha _{1}+m-1}+x_{k_{1}}^{\alpha _{1}-m+1} \\ 
x_{k_{2}}^{\alpha _{2}+1} & x_{k_{2}}^{\alpha _{2}}+x_{k_{2}}^{\alpha _{2}-2}
& \cdot \cdot \cdot \cdot \cdot \cdot \cdot \cdot \cdot \cdot & 
x_{k_{2}}^{\alpha _{2}+m-2}+x_{k_{2}}^{\alpha _{2}-m} \\ 
\cdot & \cdot & \cdot \cdot \cdot \cdot \cdot \cdot \cdot \cdot \cdot \cdot
& \cdot \\ 
\cdot & \cdot & \cdot \cdot \cdot \cdot \cdot \cdot \cdot \cdot \cdot \cdot
& \cdot \\ 
x_{k_{m}}^{\alpha _{m}-m+1} & x_{k_{m}}^{\alpha _{n}-m+2}+x_{k_{m}}^{\alpha
_{m}-m} & \cdot \cdot \cdot \cdot \cdot \cdot \cdot \cdot \cdot \cdot & 
x_{k_{m}}^{\alpha _{m}}+x_{k_{m}}^{\alpha _{m}-2m+2}
\end{array}
\right)
\end{equation*}
Set 
\begin{equation*}
\delta _{K}=\prod_{1\leq i<j\leq m}(1-\frac{x_{k_{i}}}{x_{k_{j}}}%
)\prod_{1\leq r<s\leq m}(1-\frac{1}{x_{k_{r}}x_{k_{s}}})
\end{equation*}
Then $\delta _{K}(\alpha )$ and $\delta _{K}$ belong to $\mathcal{L}_{K}.$
From a simple computation we derive the equality: 
\begin{equation}
\delta _{K}(\alpha )=\delta _{K}\cdot x_{k_{1}}^{\alpha _{1}}\cdot \cdot
\cdot x_{k_{m}}^{\alpha _{m}}.  \label{delda_factor}
\end{equation}

\noindent Consider $\eta =(\eta _{1},...,\eta _{r})$ a $r$-tuple of positive
integers summing $n.$ We define from $\eta $ the intervals $K_{1},...,K_{r}$
of $\{1,...,n\}$ by setting $K_{1}=[1,...,\eta _{1}]$ and for any $%
p=2,...,r, $ $K_{p}=[\eta _{1}+\cdot \cdot \cdot +\eta _{p-1}+1,...,\eta
_{1}+\cdot \cdot \cdot +\eta _{p}].$ Write $\delta _{\eta }=\delta
_{K_{1}}\cdot \cdot \cdot \delta _{K_{r}}.$

\noindent Set 
\begin{equation*}
\delta =\prod_{1\leq i<j\leq n}(1-\frac{x_{i}}{x_{j}})\prod_{1\leq r<s\leq
n}(1-\frac{1}{x_{r}x_{s}}).
\end{equation*}
Giving $\beta =(\beta _{1},...,\beta _{n})\in \mathbb{Z}^{n},$ set $\beta
^{(1)}=(\beta _{1},...,\beta _{\eta _{1}})$ and $\beta ^{(p)}=(\beta _{\eta
_{1}+\cdot \cdot \cdot +\eta _{p-1}+1},...,\beta _{\eta _{1}+\cdot \cdot
\cdot +\eta _{p}})$ for any $p=2,...,r.\;$The number of decompositions 
\begin{equation*}
\beta =\sum_{1\leq i<j\leq n}a_{i,j}(\varepsilon _{i}-\varepsilon
_{j})-\sum_{1\leq r<s\leq n}b_{r,s}(\varepsilon _{r}+\varepsilon _{s})
\end{equation*}
with $a_{i,j}$ and $b_{r,s}$ some positive integers is finite. Thus $\delta
^{-1}$ is well defined and belongs to $\mathcal{L}_{n}$. We introduce the
linear maps 
\begin{gather}
\left\{ 
\begin{tabular}{c}
$\mathrm{U}_{G,\eta }:\mathcal{L}_{n}\rightarrow \mathcal{F}^{G}$ \\ 
$x^{\beta }\mapsto u_{\beta ^{(1)}}^{G}\cdot \cdot \cdot u_{\beta
^{(r)}}^{G} $%
\end{tabular}
\right. ,\left\{ 
\begin{tabular}{c}
$\mathrm{U}_{G}:\mathcal{L}_{n}\rightarrow \mathcal{F}^{G}$ \\ 
$x^{\beta }\mapsto u_{\beta }^{G}$%
\end{tabular}
\right. ,\text{ }\left\{ 
\begin{tabular}{c}
$\mathrm{H}_{G}:\mathcal{L}_{n}\rightarrow \mathcal{F}^{G}$ \\ 
$x^{\beta }\mapsto h_{\beta }^{G}$%
\end{tabular}
\right.  \label{def_UH} \\
\left\{ 
\begin{tabular}{c}
$\mathrm{\Delta }_{\eta }:\mathcal{L}_{n}\rightarrow \mathcal{L}_{n}$ \\ 
$x^{\beta }\mapsto \delta _{\eta }\cdot x^{\beta }$%
\end{tabular}
\right. \text{ and }\left\{ 
\begin{tabular}{c}
$\mathrm{\nabla }:\mathcal{L}_{n}\rightarrow \mathcal{L}_{n}$ \\ 
$x^{\beta }\mapsto \delta ^{-1}\cdot x^{\beta }$%
\end{tabular}
\right.  \notag
\end{gather}
Note that these maps are not ring homomorphisms.

\begin{lemma}
\label{lem_HG}Let $K=[k_{1},...,k_{m}]\subset \{1,...,n\}$ be the interval
containing the $m$ consecutive integers $k_{1}<\cdot \cdot \cdot <k_{m}.$\
Then for any $\alpha \in \mathbb{Z}^{m}$ we have $\mathrm{H}_{G}(\delta
_{K}x_{k_{1}}^{\alpha _{1}}\cdot \cdot \cdot x_{k_{m}}^{\alpha
_{m}})=u_{\alpha }^{G}.$
\end{lemma}

\begin{proof}
To simplify the notation we set $x_{k_{i}}=y_{i}$ for any $k=1,...,m.\;$The
linear map $\mathrm{H}_{G}$ is not a ring homomorphism.\ Nevertheless, if $%
P_{1},...,P_{k}$ are polynomials respectively in the indeterminates $%
y_{1},...,y_{m}$, we have 
\begin{equation*}
\mathrm{H}_{G}(P_{1}(y_{1})\cdot \cdot \cdot P_{k}(y_{m}))=\mathrm{H}%
_{G}(P_{1}(y_{1}))\cdot \cdot \cdot \mathrm{H}_{G}(P_{k}(y_{m}))
\end{equation*}
by linearity of $\mathrm{H}_{G}.\;$We can write 
\begin{equation*}
\delta _{K}(\alpha )=\sum_{\sigma \in \mathcal{S}_{m}}(-1)^{\ell (\sigma
)}y_{\sigma (1)}^{\alpha _{1}-\sigma (1)+1}(y_{\sigma (2)}^{\alpha
_{2}-\sigma (2)+2}+y_{\sigma (2)}^{\alpha _{2}-\sigma (2)})\cdot \cdot \cdot
(y_{\sigma (m)}^{\alpha _{m}-\sigma (m)+m}+y_{\sigma (m)}^{\alpha
_{m}-\sigma (m)-m+2})
\end{equation*}
and by the previous argument 
\begin{equation*}
\mathrm{H}_{G}(\delta _{K}(\alpha ))=\sum_{\sigma \in \mathcal{S}%
_{m}}(-1)^{\ell (\sigma )}h_{\alpha _{1}-\sigma (1)+1}\cdot \cdot \cdot
(h_{\alpha _{m}-\sigma (m)+m}+h_{\alpha _{m}-\sigma (m)-m+2})=u_{\alpha }^{G}
\end{equation*}
where the last equality follows from (\ref{u_h}).\ By (\ref{delda_factor})
we have $\delta _{K}(\alpha )=\delta _{K}y^{\alpha }.\;$Thus by applying $%
\mathrm{H}_{G}$ to this equality we obtain $\mathrm{H}_{G}(\delta
_{K}y^{\alpha })=u_{\alpha }^{G}.$
\end{proof}

\begin{proposition}
\label{prop_nabla}We have 
\begin{equation*}
\mathrm{(i):U}_{G,\eta }=\mathrm{H}_{G}\circ \mathrm{\Delta }_{\eta }\text{, 
}\mathrm{(ii):H}_{G}=\mathrm{U}_{G}\circ \mathrm{\nabla }\text{ and }\mathrm{%
(iii):\mathrm{U}}_{G,\eta }=\mathrm{U}_{G}\circ \mathrm{\nabla \circ \Delta }%
_{\eta }
\end{equation*}
\end{proposition}

\begin{proof}
$\mathrm{(i):}$ Consider $\beta \in \mathbb{Z}^{n}.$ If $P_{1},...,P_{k}$
are polynomials in the indeterminates belonging respectively to the sets $%
\{x_{i_{1}}\mid i_{1}\in K_{1}\},...,\{x_{i_{r}}\mid i_{r}\in K_{r}\}$ we
have as in the proof of the previous lemma 
\begin{equation*}
\mathrm{H}_{G}(P_{1}\cdot \cdot \cdot P_{k})=\mathrm{H}_{G}(P_{1})\cdot
\cdot \cdot \mathrm{H}_{G}(P_{k}).
\end{equation*}
Since $\mathrm{\Delta }_{\eta }(x^{\beta })=\delta _{K_{1}}x^{\beta
^{(1)}}\cdot \cdot \cdot \delta _{K_{r}}x^{\beta ^{(r)}}$ where the
polynomials $\delta _{K_{p}}x^{\beta ^{(p)}},p=1,...,r$ are respectively in
the variables $\{x_{i_{1}}\mid i_{1}\in K_{1}\},...,\{x_{i_{r}}\mid i_{r}\in
K_{r}\}$ we can write 
\begin{equation*}
\mathrm{H}_{G}\circ \mathrm{\Delta }_{\eta }(x^{\beta })=\mathrm{H}%
_{G}(\delta _{K_{1}}x^{\beta ^{(1)}})\cdot \cdot \cdot \mathrm{H}_{G}(\delta
_{K_{r}}x^{\beta ^{(r)}}).
\end{equation*}
By applying Lemma \ref{lem_HG}, we derive $\mathrm{H}_{G}\circ \mathrm{%
\Delta }_{\eta }(x^{\beta })=u_{\beta ^{(1)}}^{G}\cdot \cdot \cdot u_{\beta
^{(r)}}^{G}=\mathrm{U}_{G,\eta }(x^{\beta }).$

\noindent $\mathrm{(ii):}$ Consider the linear map 
\begin{equation*}
\left\{ 
\begin{tabular}{c}
$\mathrm{\Delta }:\mathcal{L}_{n}\rightarrow \mathcal{L}_{n}$ \\ 
$x^{\beta }\mapsto \delta \cdot x^{\beta }$%
\end{tabular}
\right. .
\end{equation*}
Then by Lemma \ref{lem_HG} applied with $K=\{1,...,n\}$ we will have $%
\mathrm{H}_{G}\circ \mathrm{\Delta }=\mathrm{U}_{G}.$ Now for any $\beta \in 
\mathbb{Z}^{n}$ it is clear that $\mathrm{\Delta }\circ \mathrm{\nabla }%
(x^{\beta })=x^{\beta }.\;$This implies that $\mathrm{U}_{G}\circ \mathrm{%
\nabla }(x^{\beta })=\mathrm{H}_{G}\circ \mathrm{\Delta }\circ \mathrm{%
\nabla }(x^{\beta })=\mathrm{H}_{G}(x^{\beta })$ for any $\beta \in \mathbb{Z%
}^{n}$ and $\mathrm{(ii)}$ is proved.

\noindent $\mathrm{(iii)}$ is a straightforward consequence of $\mathrm{(i)}$
and $\mathrm{(ii)}$.
\end{proof}

\bigskip

\noindent Now we have the equality 
\begin{equation}
\delta ^{-1}\cdot \delta _{\eta }=\prod_{(i,j)\in E_{\eta }}(1-\frac{x_{i}}{%
x_{j}})^{-1}\prod_{(r,s)\in E_{\eta }}(1-\frac{1}{x_{r}x_{s}})^{-1}
\label{def_delta}
\end{equation}
where $E_{\eta }=\cup _{2\leq p\leq r}\{(i,j)\mid 1\leq i\leq \eta
_{1}+\cdot \cdot \cdot +\eta _{p-1}<j\leq n\}.$ In particular $\delta
^{-1}\cdot \delta _{\eta }$ belongs to $\mathcal{L}_{n}.$ Set 
\begin{equation}
\delta ^{-1}\cdot \delta _{\eta }=\sum_{\beta \in \mathbb{Z}^{n}}\mathcal{Q}%
^{\eta }(\beta )x^{\beta }.  \label{delta}
\end{equation}
Let $\mathbf{\mu }=(\mu ^{(1)},...,\mu ^{(r)})$ be a $r$-tuple of partitions
such that $\mu ^{(k)}$ belongs to $\mathcal{P}_{\eta _{k}}$ for any $%
k=1,...,r.$ For $G=SO_{2n+1},Sp_{2n}$ or $SO_{2n},$ the
Littlewood-Richardson coefficients $d_{\mu ^{(1)},...,\mu ^{(r)}}^{\lambda }$
are defined by the equality 
\begin{equation*}
s_{\mu ^{(1)}}^{G}\cdot \cdot \cdot s_{\mu ^{(r)}}^{G}=\sum_{\lambda \in 
\mathcal{P}_{n}}d_{\mu ^{(1)},...,\mu ^{(r)}}^{\lambda }s_{\lambda }^{G}.
\end{equation*}
This means that the coefficient $d_{\mu ^{(1)},...,\mu ^{(r)}}^{\lambda }$
gives the multiplicity of the irreducible $G$-module $V^{G}(\lambda )$ in
the tensor product $V^{G}(\mu ^{(1)})\otimes \cdot \cdot \cdot \otimes
V^{G}(\mu ^{(r)})$ thus is a nonnegative integer.

\begin{proposition}
\label{pro_d}With the above notation we have 
\begin{equation*}
d_{\mu ^{(1)},...,\mu ^{(r)}}^{\lambda }=\sum_{\sigma \in S_{n}}(-1)^{\ell
(\sigma )}\mathcal{Q}^{\eta }(\sigma \circ \lambda -\mu )
\end{equation*}
where $\mu \in \mathbb{N}^{n}$ is obtained by reading successively the parts
of the partitions $\mu ^{(1)},...,\mu ^{(r)}$ defining $\mathbf{\mu }$ from
left to right.
\end{proposition}

\begin{proof}
By $\mathrm{(iii)}$ of Proposition \ref{prop_nabla}, we have $\mathrm{%
\mathrm{U}}_{G,\eta }=\mathrm{U}_{G}\circ \mathrm{\nabla \circ \Delta }%
_{\eta }.$ Since $\mathrm{\nabla \circ \Delta }_{\eta }(x^{\mu })=\delta
^{-1}\cdot \delta _{\eta }\cdot x^{\mu },$ we obtain by (\ref{delta}) 
\begin{equation*}
\mathrm{\mathrm{U}}_{G,\eta }(x^{\mu })=u_{\mu ^{(1)}}^{G}\cdot \cdot \cdot
u_{\mu ^{(r)}}^{G}=\sum_{\beta \in \mathbb{Z}^{n}}\mathcal{Q}^{\eta }(\beta
)u_{\beta +\mu }^{G}.
\end{equation*}
Now by using Lemma \ref{straight} we derive 
\begin{equation*}
s_{\mu ^{(1)}}^{G}\cdot \cdot \cdot s_{\mu ^{(r)}}^{G}=\sum_{\lambda \in 
\mathcal{P}_{n}}\sum_{\sigma \in S_{n}}(-1)^{\ell (\sigma )}\mathcal{Q}%
^{\eta }(\sigma \circ \lambda -\mu )s_{\lambda }^{G}
\end{equation*}
and the proposition is proved.
\end{proof}

\bigskip

\noindent \textbf{Remark: }When $\mathbf{\mu }=(\mu ^{(1)},...,\mu ^{(r)})$
is a $r$-tuple of partitions such that $\mu \in \mathbb{N}^{n},$ we recover
from the above proposition that the coefficients $d_{\mu ^{(1)},...,\mu
^{(r)}}^{\lambda }$ do not depend on the Lie group $G=SO_{2n+1},Sp_{2n}$ or $%
SO_{2n}$ considered.

\subsection{A duality \label{subsec_dual} for the coefficients $d_{\protect%
\mu ^{(1)},...,\protect\mu ^{(r)}}^{\protect\lambda }$}

We define the involution $\iota $ on $\mathbb{Z}^{n}$ by setting $\iota
(\beta _{1},...,\beta _{n})=(-\beta _{n},...,-\beta _{1})$ for any $\beta
=(\beta _{1},...,\beta _{n})\in \mathbb{Z}^{n}.$ Let $\eta =(\eta
_{1},...,\eta _{r})$ be a $r$-tuple of positive integers summing $n.$ Set $%
\overline{\eta }=(\eta _{r},...,\eta _{1}).$

\begin{lemma}
\label{lem_invol_I}For any $\beta =(\beta _{1},...,\beta _{n})\in \mathbb{Z}%
^{n}$ we have 
\begin{equation*}
\mathcal{Q}^{\eta }(\beta )=\mathcal{P}^{\overline{\eta }}(\iota (\beta ))
\end{equation*}
where $\mathcal{P}^{\overline{\eta }}$ and $\mathcal{Q}^{\eta }$ are
respectively the partition functions defined in (\ref{par_So}) and (\ref
{delta}).
\end{lemma}

\begin{proof}
By abuse of notation we also denote by $\iota $ the ring automorphism of $%
\mathcal{L}_{n}$ defined by $\iota (x^{\beta })=x^{\iota (\beta )}.$ By
applying $\iota $ to the identity 
\begin{equation*}
\prod_{(i,j)\in E_{\eta }}(1-\frac{x_{i}}{x_{j}})^{-1}\prod_{(r,s)\in
E_{\eta }}(1-\frac{1}{x_{r}x_{s}})^{-1}=\sum_{\beta \in \mathbb{Z}^{n}}%
\mathcal{Q}^{\eta }(\beta )x^{\beta }
\end{equation*}
we obtain 
\begin{equation*}
\prod_{(i,j)\in E_{\overline{\eta }}}(1-\frac{x_{i}}{x_{j}}%
)^{-1}\prod_{(r,s)\in E_{\overline{\eta }}}(1-x_{r}x_{s})^{-1}=\sum_{\beta
\in \mathbb{Z}^{n}}\mathcal{Q}^{\eta }(\beta )x^{\iota (\beta )}=\sum_{\beta
\in \mathbb{Z}^{n}}\mathcal{P}^{\overline{\eta }}(\beta )x^{\beta }
\end{equation*}
where $E_{\overline{\eta }}=\cup _{2\leq p\leq r}\{(i,j)\mid 1\leq i\leq 
\overline{\eta }_{1}+\cdot \cdot \cdot +\overline{\eta }_{p-1}<j\leq n\}.$
This implies $\mathcal{Q}^{\eta }(\beta )=\mathcal{P}^{\overline{\eta }%
}(\iota (\beta ))$ for any $\beta \in \mathbb{Z}^{n}.$
\end{proof}

\bigskip\ \ 

\noindent Given $\sigma \in S_{n},$ denote by $\overline{\sigma }$ the
permutation defined by 
\begin{equation*}
\overline{\sigma }(k)=\sigma (n-k+1).
\end{equation*}
For any $i\in \{1,...,n-1\},$ we have $\overline{s}_{i}=s_{n-i}.\;$The
following Lemma is straightforward:

\begin{lemma}
\label{lem_invbar}The map $\sigma \rightarrow \overline{\sigma }$ is an
involution of the group $S_{n}$.\ Moreover we have $\sigma (\iota (\beta
))=\iota (\overline{\sigma }(\beta ))$ and $\ell (\sigma )=\ell (\overline{%
\sigma })$ for any $\beta \in \mathbb{Z}^{n},\sigma \in S_{n}.$
\end{lemma}

\begin{lemma}
Let $\lambda ,\mu $ two partitions of length $n$ and $\sigma \in S_{n}.$
Then 
\begin{equation*}
(-1)^{\ell (\sigma )}\mathcal{Q}^{\eta }(\sigma (\lambda +\rho )-(\mu +\rho
))=(-1)^{\ell (\overline{\sigma })}\mathcal{P}^{\overline{\eta }}(\overline{%
\sigma }(\iota (\lambda )+\rho )-(\iota (\mu )+\rho ))
\end{equation*}
\end{lemma}

\begin{proof}
It suffices to prove the identity 
\begin{equation*}
\mathcal{Q}^{\eta }(\sigma (\lambda +\rho )-(\mu +\rho ))=\mathcal{P}^{%
\overline{\eta }}(\overline{\sigma }(\iota (\lambda )+\rho )-(\iota (\mu
)+\rho ))
\end{equation*}
Set $P=\mathcal{P}^{\overline{\eta }}(\overline{\sigma }(\iota (\lambda
)+\rho )-(\iota (\mu )+\rho )).\;$From the above Lemma we deduce 
\begin{equation*}
P=\mathcal{P}^{\overline{\eta }}(\iota (\sigma (\lambda ))+\overline{\sigma }%
(\rho )-\iota (\mu )-\rho ).
\end{equation*}
Now an immediate computation shows that $\overline{\sigma }(\rho )-\rho
=\iota (\sigma (\rho )-\rho ).\;$Thus we derive 
\begin{equation*}
P=\mathcal{P}^{\overline{\eta }}(\iota (\sigma (\lambda +\rho )-\mu -\rho ))=%
\mathcal{Q}^{\eta }(\sigma (\lambda +\rho )-\mu -\rho )
\end{equation*}
where the last equality follows from Lemma \ref{lem_invol_I}.
\end{proof}

\bigskip

\noindent Consider $\lambda \in \mathcal{P}_{n}$ and $\mathbf{\mu }=(\mu
^{(1)},...,\mu ^{(r)})$ a $r$-tuple of partitions such that $\mu ^{(k)}$
belongs to $\mathcal{P}_{\eta _{k}}$ for any $k=1,...,r.$ Recall that $\mu
\in \mathbb{N}^{n}$ is the $n$-tuple obtained by reading successively the
parts of the partitions $\mu ^{(1)},...,\mu ^{(r)}$ defining $\mathbf{\mu }$
from left to right. Let $a$ be the minimal integer such that 
\begin{equation}
\widehat{\lambda }=(a-\lambda _{n},...,a-\lambda _{1})\text{ and }\widehat{%
\mu }=(a-\mu _{n},...,a-\mu _{1})  \label{def-lambdahat}
\end{equation}
belong to $\mathbb{N}^{n}.\;$Then $\widehat{\lambda }$ is a partition of
length $n$.\ Set $\overline{\eta }=(\overline{\eta }_{1},...,\overline{\eta }%
_{r})$ and denote by $\widehat{\mathbf{\mu }}=(\widehat{\mu }^{(1)},...,%
\widehat{\mu }^{(r)})$ the $r$-tuple of partitions such that $\widehat{\mu }%
^{(1)}=(\mu _{1},...,\mu _{\overline{\eta }_{1}})\in \mathcal{P}_{\overline{%
\eta }_{1}}$ and $\widehat{\mu }^{(p)}=(\mu _{\overline{\eta }_{1}+\cdot
\cdot \cdot +\overline{\eta }_{p}+1},...,\mu _{\overline{\eta }_{1}+\cdot
\cdot \cdot +\overline{\eta }_{p}})\in \mathcal{P}_{\overline{\eta }_{p}}$
for any $k=2,...,r.$ The following proposition shows that the
Littlewood-Richardson coefficients $d_{\mu ^{(1)},...,\mu ^{(r)}}^{\lambda }$
defined above are branching coefficients associated to the restriction from
the orthogonal group $SO_{2n}$ to the subgroup $SO_{\overline{\eta }}\simeq
SO_{\overline{\eta }_{1}}\times \cdot \cdot \cdot SO_{\overline{\eta }_{r}}$
defined in \ref{subsec_lit}.

\begin{proposition}
\label{th_dual1}With the above notation, we have for any integer $k\geq 
\frac{\left| \mu \right| -\left| \lambda \right| }{2}$%
\begin{equation}
d_{\mu ^{(1)},...,\mu ^{(r)}}^{\lambda }=[V(\widehat{\lambda }+k\kappa
)^{SO_{2n}}:V(\widehat{\mu }+k\kappa )^{SO_{\overline{\eta }}}].  \label{th}
\end{equation}
\end{proposition}

\begin{proof}
It follows from the definition of $d_{\mu ^{(1)},...,\mu ^{(r)}}^{\lambda }$
and the above lemma that 
\begin{equation*}
d_{\mu ^{(1)},...,\mu ^{(r)}}^{\lambda }=\sum_{\sigma \in S_{n}}(-1)^{\ell
(\sigma )}\mathcal{Q}^{\eta }(\sigma (\lambda +\rho )-\mu -\rho )=\sum_{%
\overline{\sigma }\in \mathcal{S}_{n}}(-1)^{\ell (\overline{\sigma })}%
\mathcal{P}^{\overline{\eta }}(\overline{\sigma }(\iota (\lambda )+\rho
))-(\iota (\mu )+\rho )).
\end{equation*}
Then by Lemma \ref{lem_invbar} we obtain 
\begin{equation*}
d_{\mu ^{(1)},...,\mu ^{(r)}}^{\lambda }=\sum_{\sigma \in \mathcal{S}%
_{n}}(-1)^{\ell (\sigma )}\mathcal{P}^{\overline{\eta }}(\sigma (\iota
(\lambda )+\rho ))-(\iota (\mu )+\rho )).
\end{equation*}
We have $\sigma (\iota (\lambda )+\rho +a\kappa )=\sigma (\iota (\lambda
)+\rho )+a\kappa $ since $\sigma \in \mathcal{S}_{n}$ and $\kappa
=(1,...,1). $ So we can write 
\begin{equation*}
d_{\mu ^{(1)},...,\mu ^{(r)}}^{\lambda }=\sum_{\sigma \in \mathcal{S}%
_{n}}(-1)^{\ell (\sigma )}\mathcal{P}^{\overline{\eta }}(\sigma (\iota
(\lambda )+a\kappa +\rho ))-(\iota (\mu )+a\kappa +\rho )).
\end{equation*}
Since $\widehat{\lambda }=\iota (\lambda )+a\kappa $ and $\widehat{\mu }%
=\iota (\mu )+a\kappa $ we derive 
\begin{equation*}
d_{\mu ^{(1)},...,\mu ^{(r)}}^{\lambda }=\sum_{\sigma \in \mathcal{S}%
_{n}}(-1)^{\ell (\sigma )}\mathcal{P}^{\overline{\eta }}(\sigma (\widehat{%
\lambda }+\rho )-(\widehat{\mu }+\rho )).
\end{equation*}
Now by using Lemma \ref{lem_tech} with $\mathcal{M}=\mathcal{P}^{\overline{%
\eta }}$ and Proposition \ref{prop_rest_diret_So} we obtain 
\begin{equation*}
\sum_{\sigma \in \mathcal{S}_{n}}(-1)^{\ell (\sigma )}\mathcal{P}^{\overline{%
\eta }}(\sigma (\widehat{\lambda }+\rho )-(\widehat{\mu }+\rho ))=[V(%
\widehat{\lambda }+k\kappa )^{SO_{2n}}:V(\widehat{\mu }+k\kappa )^{SO_{%
\overline{\eta }}}]
\end{equation*}
for any integer $k\geq \frac{\left| \widehat{\lambda }\right| -\left| 
\widehat{\mu }\right| }{2}=\frac{\left| \mu \right| -\left| \lambda \right| 
}{2}$ and the Proposition is proved.
\end{proof}

\bigskip

\noindent \textbf{Remarks:}

\noindent $\mathrm{(i):}$ The previous proposition can be regarded as an
analogue for the Lie groups $SO_{2n+1},Sp_{2n}$ and $SO_{2n}$ of the duality 
\begin{equation}
c_{\mu ^{(1)},...,\mu ^{(r)}}^{\lambda }=[V(\lambda )^{GL_{n}}:V(\mu
)^{GL_{\eta }}]  \label{dualA}
\end{equation}
between the Littlewood-Richardson coefficient $c_{\mu ^{(1)},...,\mu
^{(r)}}^{\lambda }$ giving the multiplicity of $V(\lambda )^{GL_{n}}$ in the
tensor product $V^{GL_{n}}(\mu ^{(1)})\otimes \cdot \cdot \cdot \otimes
V^{GL_{n}}(\mu ^{(r)})$ and the branching coefficient of the restriction of $%
V(\lambda )^{GL_{n}}$ to the Levi subgroup $GL_{\eta }=GL_{\eta _{1}}\times
\cdot \cdot \cdot \times GL_{\eta _{k}}.$ Note that this duality can be
proved by similar methods than those used in this paragraph starting from
Jacobi-Trudi's determinantal expression for the Schur function $s_{\lambda
}^{GL_{n}}=\mathrm{char}(V^{GL_{n}}(\lambda ))$ instead of (\ref{u_h}).

\noindent $\mathrm{(ii):}$ When all the $\mu ^{(k)}$'s are row partitions, (%
\ref{dualA}) simply express the Schur-Weyl duality between the dimension of
the weight space $\mu $ in $V^{GL_{n}}(\lambda )$ and the multiplicity of $%
V^{GL_{n}}(\lambda )$ in the tensor product of the symmetric powers of the
vector representation of $GL_{n}$ associated to $\mu .$ Similarly, in this
particular case, (\ref{th}) reduce to the duality already observed in \cite
{lec}.

\subsection{Quantization of the coefficients $d_{\protect\mu ^{(1)},...,%
\protect\mu ^{(r)}}^{\protect\lambda }$}

Consider $\lambda \in \mathcal{P}_{n}$, $\eta =(\eta _{1},...,\eta _{r})$ a $%
r$-tuple of positive integers summing $n$ and $\mathbf{\mu }=(\mu
^{(1)},...,\mu ^{(r)})$ a $r$-tuple of partitions such that $\mu ^{(k)}$
belongs to $\mathcal{P}_{\eta _{k}}$ for any $k=1,...,r.$ We have seen in
paragraph \ref{subsec_quant} that it is possible to define natural $q$%
-analogues of the Littlewood-Richardson coefficients $c_{\mu ^{(1)},...,\mu
^{(r)}}^{\lambda }$ from the duality (\ref{dualA}) by setting 
\begin{equation}
c_{\mu ^{(1)},...,\mu ^{(r)}}^{\lambda }(q)=K_{\lambda ,\mu }^{GL_{n},I}(q).
\label{qdualA}
\end{equation}
where the polynomials $K_{\lambda ,\mu }^{GL_{n},I}(q)$ are the generalized
Lusztig $q$-analogues of Definition \ref{def_K}. When the $\mu ^{(k)}$'s are
rectangular partitions of decreasing heights, the polynomials $c_{\mu
^{(1)},...,\mu ^{(r)}}^{\lambda }(q)$ have nonnegative coefficients by
Theorem\ref{TH_broer}. It is tempting to define $q$-analogues of the
coefficients $d_{\mu ^{(1)},...,\mu ^{(r)}}^{\lambda }$ by setting 
\begin{equation*}
d_{\mu ^{(1)},...,\mu ^{(r)}}^{\lambda }(q)=\sum_{\sigma \in
S_{n}}(-1)^{\ell (\sigma )}\mathcal{Q}_{q}^{\eta }(\sigma \circ \lambda -\mu
)
\end{equation*}
where the $q$-partition function $\mathcal{Q}_{q}^{\eta }$ verifies 
\begin{equation*}
\prod_{(i,j)\in E_{\eta }}(1-q\frac{x_{i}}{x_{j}})^{-1}\prod_{(r,s)\in
E_{\eta }}(1-q\frac{1}{x_{r}x_{s}})^{-1}=\sum_{\beta \in \mathbb{Z}^{n}}%
\mathcal{Q}_{q}^{\eta }(\beta )x^{\beta }.
\end{equation*}
Unfortunately the polynomials $d_{\mu ^{(1)},...,\mu ^{(r)}}^{\lambda }(q)$
have not nonnegative coefficients in general. This is in particular the case
for $\lambda =(1,1,1,0,0),$ $\mu ^{(1)}=(5),$ $\mu ^{(2)}=(4,4)$ and $\mu
^{(3)}=(2,2)$ where we have $d_{\mu ^{(1)},\mu ^{(2)},\mu ^{(3)}}^{\lambda
}(q)=q^{11}-q^{8}.$

\subsection{The $q$-analogues $\frak{D}_{\protect\mu ^{(1)},...,\protect\mu
^{(r)}}^{\protect\lambda ,G}(q)$}

\noindent For $G=SO_{2n+1},Sp_{2n}$ or $SO_{2n}$, there exist coefficients
associated to the decomposition of a tensor product of $G$-modules into its
irreducible components admitting a natural quantization with nonnegative
coefficients. For any partition $\nu \in \mathcal{P}_{n},$ write $\frak{V}%
^{G}(\nu )$ for the restriction from the irreducible $GL_{N}$-module of
highest weight $\nu $ to $G.$ Consider $\eta =(\eta _{1},...,\eta _{r})$ a $%
r $-tuple of positive integers summing $n$ and $\mathbf{\mu }=(\mu
^{(1)},...,\mu ^{(r)})$ a $r$-tuple of partitions such that $\mu ^{(k)}$
belongs to $\mathcal{P}_{\eta _{k}}$ for any $k=1,...,r.$ Given $\lambda \in 
\mathcal{P}_{n},$ the coefficients $\frak{D}_{\mu ^{(1)},...,\mu
^{(r)}}^{\lambda ,G}$ are defined as the multiplicities of $V^{G}(\lambda )$
in $\frak{V}^{G}(\mu ^{(1)})\otimes \cdot \cdot \cdot \otimes \frak{V}%
^{G}(\mu ^{(r)}),$ that is we have 
\begin{equation*}
\frak{V}^{G}(\mu ^{(1)})\otimes \cdot \cdot \cdot \otimes \frak{V}^{G}(\mu
^{(r)})\simeq \bigoplus_{\lambda \in \mathcal{P}_{n}}V^{G}(\lambda )^{\oplus 
\frak{D}_{\mu ^{(1)},...,\mu ^{(r)}}^{\lambda ,G}}.
\end{equation*}
Contrary to the coefficients $d_{\mu ^{(1)},...,\mu ^{(r)}}^{\lambda },$ the
coefficients $\frak{D}_{\mu ^{(1)},...,\mu ^{(r)}}^{\lambda ,G}$ depend on
the Lie group $G$ considered.\ Set 
\begin{equation*}
\Omega ^{G}=\left\{ 
\begin{tabular}{l}
$\prod_{1\leq r<s\leq n}(1-\frac{1}{x_{r}x_{s}})$ for $G=Sp_{2n}$ \\ 
$\prod_{1\leq r\leq s\leq n}(1-\frac{1}{x_{r}x_{s}})$ for $G=SO_{2n}$ \\ 
$\prod_{1\leq r<s\leq n}(1-\frac{1}{x_{r}x_{s}})\prod_{1\leq i\leq n}(1-%
\frac{1}{xi})$ for $G=SO_{2n+1}$%
\end{tabular}
\right. .
\end{equation*}
Denote by $\frak{Q}^{\eta ,G}$ the partition function defined for $%
G=Sp_{2n},SO_{2n}$ and $SO_{2n+1}$ by the identities 
\begin{equation}
\prod_{(i,j)\in E_{\eta }}(1-\frac{x_{i}}{x_{j}})^{-1}\Omega
_{G}=\sum_{\beta \in \mathbb{Z}^{n}}\frak{Q}^{\eta ,G}(\beta )x^{\beta }
\label{dec_omega}
\end{equation}
where $E_{\eta }=\cup _{2\leq p\leq r}\{(i,j)\mid 1\leq i\leq \eta
_{1}+\cdot \cdot \cdot +\eta _{p-1}<j\leq n\}.$

\begin{proposition}
With the above notation we have 
\begin{equation*}
\frak{D}_{\mu ^{(1)},...,\mu ^{(r)}}^{\lambda ,G}=\sum_{\sigma \in
S_{n}}(-1)^{\ell (\sigma )}\frak{Q}^{\eta ,G}(\sigma \circ \lambda -\mu )
\end{equation*}
\end{proposition}

\begin{proof}
Consider $m\in \mathbb{N}$ and for any $\alpha =(\alpha _{1},...,\alpha
_{m}) $ set 
\begin{equation*}
v_{\alpha }=\det \left( 
\begin{array}{cccc}
h_{\alpha _{1}}^{GL_{n}} & h_{\alpha _{1}+1}^{GL_{n}} & \cdot \cdot \cdot
\cdot \cdot \cdot \cdot \cdot \cdot \cdot & h_{\alpha _{1}+m-1}^{GL_{n}} \\ 
h_{\alpha _{2}-1}^{GL_{n}} & h_{\alpha _{2}}^{GL_{n}} & \cdot \cdot \cdot
\cdot \cdot \cdot \cdot \cdot \cdot \cdot & h_{\alpha _{2}+m-2}^{GL_{n}} \\ 
\cdot & \cdot & \cdot \cdot \cdot \cdot \cdot \cdot \cdot \cdot \cdot \cdot
& \cdot \\ 
\cdot & \cdot & \cdot \cdot \cdot \cdot \cdot \cdot \cdot \cdot \cdot \cdot
& \cdot \\ 
h_{\alpha _{m}-m+1}^{GL_{n}} & h_{\alpha _{m}-m+2}^{GL_{n}} & \cdot \cdot
\cdot \cdot \cdot \cdot \cdot \cdot \cdot \cdot & h_{\alpha _{m}}^{GL_{n}}
\end{array}
\right) .
\end{equation*}
As in Lemma \ref{straight}, we have from Jacobi-Trudi's determinantal
expression of the Schur function $s_{\lambda }^{GL_{n}}$%
\begin{equation}
v_{\alpha }=\left\{ 
\begin{tabular}{l}
$(-1)^{\ell (\sigma )}s_{\lambda }^{GL_{n}}$ if there exists $\sigma \in 
\mathcal{S}_{n}$ and $\lambda \in \mathcal{P}_{n}$ such that $\sigma \circ
\alpha =\lambda $ \\ 
$0$ otherwise
\end{tabular}
\right. .  \label{strai2}
\end{equation}
Let $K=[k_{1},...,k_{m}]\subset \{1,...,n\}$ be the interval containing the $%
m$ consecutive integers $k_{1}<\cdot \cdot \cdot <k_{m}.$\ Consider the
Laurent polynomial $\frak{d}_{K}=\prod_{1\leq i<j\leq m}(1-\frac{x_{k_{i}}}{%
x_{k_{j}}})\;$(see paragraph \ref{subsec_Jacobi-trudi}) Then we derive as in
Lemma \ref{lem_HG} the equality $\mathrm{H}_{GL_{n}}(\frak{d}_{K}\cdot
x^{\alpha })=v_{\alpha }^{G}.\;$Now observe that for any integer $k,$ we
have $h_{k}^{SP_{2n}}=h_{k}^{GL_{2n}},$ $%
h_{k}^{SO_{2n}}=h_{k}^{GL_{2n}}-h_{k-2}^{GL_{2n}}$ and $%
h_{k}^{SO_{2n+1}}=h_{k}^{GL_{2n+1}}-h_{k-1}^{GL_{2n+1}}.$ This permits to
express the determinant (\ref{u_h}) in terms of the $h_{k}^{GL_{2n}},$ $k\in 
\mathbb{Z}$. Set 
\begin{equation*}
\frak{d}^{G}=\prod_{1\leq i<j\leq n}(1-\frac{x_{i}}{x_{j}})\Omega _{G}.
\end{equation*}
For any $\beta =(\beta _{1},...,\beta _{n})\in \mathbb{Z}^{n}$ and for each
Lie group $Sp_{2n},SO_{2n+1}$ and $SO_{2n}$ we will have 
\begin{equation*}
u_{\beta }^{G}=\mathrm{H}_{GL_{n}}(\frak{d}^{G}\cdot x^{\beta }).
\end{equation*}
Now we define the intervals $K_{1},...,K_{r}$ of $\{1,...,n\}$ from $\eta $
by setting $K_{1}=[1,...,\eta _{1}]$ and for any $p=2,...,r,$ $K_{p}=[\eta
_{1}+\cdot \cdot \cdot +\eta _{p-1}+1,...,\eta _{1}+\cdot \cdot \cdot +\eta
_{p}].$ Consider the linear maps 
\begin{equation}
\left\{ 
\begin{tabular}{c}
$\mathrm{V}_{\eta }:\mathcal{L}_{n}\rightarrow \mathcal{F}_{n}^{GL_{n}}$ \\ 
$x^{\beta }\mapsto v_{\beta ^{(1)}}\cdot \cdot \cdot v_{\beta ^{(r)}}$%
\end{tabular}
\right. \left\{ 
\begin{tabular}{c}
$\mathrm{\Phi }_{\eta }:\mathcal{L}_{n}\rightarrow \mathcal{L}_{n}$ \\ 
$x^{\beta }\mapsto \frak{d}_{\eta }\cdot x^{\beta }$%
\end{tabular}
\right. \text{ and }\left\{ 
\begin{tabular}{c}
$\mathrm{\Psi }^{G}:\mathcal{L}_{n}\rightarrow \mathcal{L}_{n}$ \\ 
$x^{\beta }\mapsto (\frak{d}^{G})^{-1}\cdot x^{\beta }$%
\end{tabular}
\right.  \notag
\end{equation}
where $\frak{d}_{\eta }=\frak{d}_{K_{1}}\cdot \cdot \cdot \frak{d}_{K_{r}}$
and $\beta ^{(1)},...,\beta ^{(r)}$ are defined as in (\ref{def_UH}).$\;$%
Then we obtain as in Proposition \ref{prop_nabla} the identities $\mathrm{%
\mathrm{V}}_{\eta }=\mathrm{H}_{GL_{n}}\circ \mathrm{\Phi }_{\eta }$, $%
\mathrm{H}_{GL_{n}}=\mathrm{U}_{G}\circ \mathrm{\Psi }^{G}$ which lead to $%
\mathrm{\mathrm{V}}_{\eta }=\mathrm{U}_{G}\circ \mathrm{\Psi }^{G}\circ 
\mathrm{\Phi }_{\eta }.$ By (\ref{dec_omega}) we have 
\begin{equation*}
(\frak{d}^{G})^{-1}\cdot \frak{d}_{\eta }=\prod_{(i,j)\in E_{\eta }}(1-\frac{%
x_{i}}{x_{j}})^{-1}\Omega ^{G}=\sum_{\beta \in \mathbb{Z}^{n}}\frak{Q}^{\eta
,G}(\beta )x^{\beta }
\end{equation*}
thus we can write $\mathrm{\mathrm{V}}_{\eta }(x^{\mu })=\mathrm{U}%
_{G}(\sum_{\beta \in \mathbb{Z}^{n}}^{\eta }\frak{Q}^{\eta ,G}(\beta
)x^{\beta +\mu })=\sum_{\beta \in \mathbb{Z}^{n}}^{\eta }\frak{Q}^{\eta
,G}(\beta )u_{\beta +\mu }^{G}.$ Then by using the straightening law (\ref
{strai2}) we derive 
\begin{equation*}
s_{\mu ^{(1)}}^{GL_{n}}\cdot \cdot \cdot s_{\mu
^{(r)}}^{GL_{n}}=\sum_{\lambda \in \mathcal{P}_{n}}\sum_{\sigma \in
S_{n}}(-1)^{\ell (\sigma )}\frak{Q}^{\eta ,G}(\sigma \circ \lambda -\mu
)s_{\lambda }^{G}
\end{equation*}
which establishes the proposition.
\end{proof}

\bigskip

\noindent The coefficients $\frak{D}_{\mu ^{(1)},...,\mu ^{(r)}}^{\lambda
,G} $ admit the natural quantization 
\begin{equation*}
\frak{D}_{\mu ^{(1)},...,\mu ^{(r)}}^{\lambda ,G}(q)=\sum_{\sigma \in
S_{n}}(-1)^{\ell (\sigma )}\frak{Q}_{q}^{\eta ,G}(\sigma \circ \lambda -\mu )
\end{equation*}
where $\frak{Q}_{q}^{\eta ,G}$ is defined for $G=Sp_{2n},SO_{2n}$ and $%
SO_{2n+1}$ respectively by the identities 
\begin{equation}
\left\{ 
\begin{tabular}{l}
$\prod_{(i,j)\in E_{\eta }}(1-q\frac{x_{i}}{x_{j}})^{-1}\prod_{1\leq r<s\leq
n}(1-q\frac{1}{x_{r}x_{s}})^{-1}=\sum_{\beta \in \mathbb{Z}^{n}}\frak{Q}%
_{q}^{\eta ,Sp_{2n}}(\beta )x^{\beta }$ \\ 
$\prod_{(i,j)\in E_{\eta }}(1-q\frac{x_{i}}{x_{j}})^{-1}\prod_{1\leq r\leq
s\leq n}(1-q\frac{1}{x_{r}x_{s}})^{-1}=\sum_{\beta \in \mathbb{Z}^{n}}\frak{Q%
}_{q}^{\eta ,SO_{2n}}(\beta )x^{\beta }$ \\ 
$\prod_{(i,j)\in E_{\eta }}(1-q\frac{x_{i}}{x_{j}})^{-1}\prod_{1\leq r<s\leq
n}(1-q\frac{1}{x_{r}x_{s}})^{-1}\prod_{1\leq i\leq n}(1-q^{2}\frac{1}{x_{i}}%
)^{-1}=\sum_{\beta \in \mathbb{Z}^{n}}\frak{Q}_{q}^{\eta ,SO_{2n+1}}(\beta
)x^{\beta }$%
\end{tabular}
\right. .  \label{def_q-det}
\end{equation}
There exists a duality between the polynomials $\frak{D}_{\mu ^{(1)},...,\mu
^{(r)}}^{\lambda ,G}(q)$ and the $q$-analogues of Theorem \ref{prop_dec_K_c}%
. Surprisingly this duality does not relate the polynomials $\frak{D}_{\mu
^{(1)},...,\mu ^{(r)}}^{\lambda ,G}(q)$ to the polynomials $\widetilde{K}_{%
\widehat{\lambda },\widehat{\mu }}^{G,I}(q)$ but to the polynomials $%
\widetilde{K}_{\widehat{\lambda },\widehat{\mu }}^{\widehat{G},I}(q)$ where 
\begin{equation*}
\widehat{G}=\left\{ 
\begin{tabular}{l}
$Sp_{2n}$ if $G=SO_{2n}$ \\ 
$SO_{2n}$ if $G=Sp_{2n}$ \\ 
$SO_{2n+1}$ if $G=SO_{2n+1}$%
\end{tabular}
\right. .
\end{equation*}
The polynomials $\frak{D}_{\mu ^{(1)},...,\mu ^{(r)}}^{\lambda ,G}(q)$ also
decompose as nonnegative combinations of Poincar\'{e} polynomials. Set $%
\overline{\eta }=(\overline{\eta }_{1},...,\overline{\eta }_{r})$.\ Then $G$
and $GL_{n}$ contain Levi subgroups $L_{G,I}$ and $L_{GL_{n},I}$ isomorphic
to $GL_{\overline{\eta }_{1}}\times \cdot \cdot \cdot \times GL_{\overline{%
\eta }_{r}}.\;$We have $I=\{0<\alpha _{i}<\overline{\eta }_{1}\}\cup _{1\leq
p\leq r-1}\{\alpha _{i}\mid \overline{\eta }_{1}+\cdot \cdot \cdot +%
\overline{\eta }_{p}<i<\overline{\eta }_{1}+\cdot \cdot \cdot +\overline{%
\eta }_{p+1}\}.$

\begin{theorem}
\label{th_qdual}With the previous notation we have 
\begin{equation}
\frak{D}_{\mu ^{(1)},...,\mu ^{(r)}}^{\lambda ,G}(q)=\widetilde{K}_{\widehat{%
\lambda },\widehat{\mu }}^{\widehat{G},I}(q)=q^{\tfrac{\left| \mu \right|
-\left| \lambda \right| }{2}}\sum_{\nu \in \mathcal{P}_{n}}[V(\nu )^{%
\widehat{G}}:V(\lambda )^{GL_{n}}]K_{\nu ,\mu }^{GL_{n},I}(q)  \label{Dfrak}
\end{equation}
where $\widehat{\lambda },\widehat{\mu }$ are defined as in Theorem \ref
{th_dual1}.
\end{theorem}

\begin{proof}
We use the notation of paragraph \ref{subsec_dual}. We only give the proof
for $G=SP_{2n}$ or $G=SO_{2n}.$ The proof is essentially the same for $%
G=SO_{2n+1}.\;$By proceeding as in the proof of Lemma \ref{lem_invol_I}, we
establish, for any $\beta =(\beta _{1},...,\beta _{n})\in \mathbb{Z}^{n}$,
the identity 
\begin{equation*}
\frak{Q}_{q}^{\eta ,G}(\beta )=\mathcal{P}_{q}^{\widehat{G},I}(\iota (\beta
)).
\end{equation*}
We have 
\begin{equation*}
\frak{D}_{\mu ^{(1)},...,\mu ^{(r)}}^{\lambda ,G}(q)=\sum_{\sigma \in 
\mathcal{S}_{n}}(-1)^{\ell (\sigma )}\frak{Q}_{q}^{\eta ,G}(\sigma (\lambda
+\rho )-(\mu +\rho ))=\sum_{\sigma \in \mathcal{S}_{n}}(-1)^{\ell (\sigma )}%
\mathcal{P}_{q}^{\widehat{G},I}(\iota (\sigma (\lambda +\rho )-(\mu +\rho
))).
\end{equation*}
Then the equality $\frak{D}_{\mu ^{(1)},...,\mu ^{(r)}}^{\lambda ,G}(q)=%
\widetilde{K}_{\widehat{\lambda },\widehat{\mu }}^{\widehat{G},I}(q)$ is
obtained by using essentially the same arguments as in the proof of theorem 
\ref{th_dual1}.

\noindent We deduce from (\ref{RG}) the expansions 
\begin{equation*}
\prod_{1\leq r<s\leq n}(1-q\frac{1}{x_{r}x_{s}})^{-1}=q^{\left| \delta
\right| /2}\sum_{\delta \in \mathbb{N}^{n}}r_{SO_{2n}}(\delta )x^{-\delta }%
\text{ and }\prod_{1\leq r\leq s\leq n}(1-q\frac{1}{x_{r}x_{s}}%
)^{-1}=q^{\left| \delta \right| /2}\sum_{\delta \in \mathbb{N}%
^{n}}r_{Sp_{2n}}(\delta )x^{-\delta }.
\end{equation*}
Consider $\beta \in \mathbb{Z}^{n}.$ By proceeding as for Lemma \ref
{lem_util} we derive from (\ref{def_q-det}) the equality 
\begin{equation*}
\frak{Q}_{q}^{\eta ,G}(\beta )=q^{-\left| \beta \right| /2}\sum_{\delta \in 
\mathbb{N}^{n},\left| \delta \right| =-\left| \beta \right| }r_{\widehat{G}%
}(\delta )\mathcal{P}_{q}^{GL_{n},I}(\beta +\delta ).
\end{equation*}
This implies the decomposition 
\begin{equation*}
\frak{D}_{\mu ^{(1)},...,\mu ^{(r)}}^{\lambda ,G}(q)=\sum_{\sigma \in 
\mathcal{S}_{n}}(-1)^{\ell (\sigma )}\sum_{\delta \in \mathbb{N}^{n},\left|
\delta \right| =-\left| \beta \right| }r_{\widehat{G}}(\delta )q^{-\left|
\beta \right| /2}\mathcal{P}_{q}^{GL_{n},I}(\sigma (\lambda +\rho +\delta
)-(\mu +\rho ))
\end{equation*}
where $\beta =\sigma (\lambda +\rho )-(\mu +\rho )$ in the second sum.\
Since $\left| \beta \right| =\left| \lambda \right| -\left| \mu \right| ,$
we have 
\begin{equation*}
\frak{D}_{\mu ^{(1)},...,\mu ^{(r)}}^{\lambda ,G}(q)=q^{\frac{\left| \mu
\right| -\left| \lambda \right| }{2}}\sum_{\sigma \in \mathcal{S}%
_{n}}(-1)^{\ell (\sigma )}\sum_{\delta \in \mathbb{N}^{n},\left| \delta
\right| =\left| \mu \right| -\left| \lambda \right| }r_{\widehat{G}}(\delta )%
\mathcal{P}_{q}^{GL_{n},I}(\sigma (\lambda +\rho +\sigma ^{-1}(\delta
))-(\mu +\rho ))
\end{equation*}
For any $\sigma \in \mathcal{S}_{n},$ we have $\sigma ^{-1}(\mathbb{N}^{n})=%
\mathbb{N}^{n}$ and $r_{\widehat{G}}(\delta )=r_{\widehat{G}}(\sigma (\delta
))$.\ Thus we obtain 
\begin{multline*}
\frak{D}_{\mu ^{(1)},...,\mu ^{(r)}}^{\lambda ,G}(q)=q^{\frac{\left| \mu
\right| -\left| \lambda \right| }{2}}\sum_{\sigma \in \mathcal{S}%
_{n}}(-1)^{\ell (\sigma )}\sum_{\delta \in \mathbb{N}^{n},\left| \delta
\right| =\left| \mu \right| -\left| \lambda \right| }r_{\widehat{G}}(\delta )%
\mathcal{P}_{q}^{GL_{n},I}(\sigma (\lambda +\rho +\delta )-(\mu +\rho ))= \\
q^{\frac{\left| \mu \right| -\left| \lambda \right| }{2}}\sum_{\delta \in 
\mathbb{N}^{n},\left| \delta \right| =\left| \mu \right| -\left| \lambda
\right| }r_{\widehat{G}}(\delta )K_{\lambda +\delta ,\mu }^{GL_{n},I}(q).
\end{multline*}
Now by Lemma \ref{strilaw}, $K_{\lambda +\delta ,\mu }^{GL_{n},I}(q)=0$ or
there exits $\sigma \in \mathcal{S}_{n}$ and $\nu \in \widetilde{\mathcal{P}}%
_{n}$ such that $\nu =\sigma ^{-1}\circ (\lambda +\delta ).$ Since $\lambda
+\delta $ has positive coordinates, $\nu $ cannot have negative coordinates
and thus $\nu \in \mathcal{P}_{n}.\;$Then we have $\left| \nu \right|
=\left| \lambda \right| +\left| \delta \right| =\left| \mu \right| $ and $%
\delta =\sigma (\nu +\rho )-\rho -\lambda $ and it follows that 
\begin{equation*}
\frak{D}_{\mu ^{(1)},...,\mu ^{(r)}}^{\lambda ,G}(q)=q^{\frac{\left| \mu
\right| -\left| \lambda \right| }{2}}\sum_{\nu \in \mathcal{P}%
_{n}}\sum_{\sigma \in \mathcal{S}_{n}}(-1)^{\ell (\sigma )}r_{\widehat{G}%
}(\sigma (\nu +\rho )-\rho -\lambda )K_{\nu ,\mu }^{GL_{n},I}(q).
\end{equation*}
Now we deduce from Lemma \ref{lem_tech} applied with $\mathcal{M}=r_{%
\widehat{G}}$ and Proposition \ref{prop_mul_sum} that the equality 
\begin{equation*}
\sum_{\sigma \in \mathcal{S}_{n}}(-1)^{\ell (\sigma )}r_{\widehat{G}}(\sigma
(\nu +\rho )-(\lambda +\rho ))=[V(\nu +k\kappa )^{\widehat{G}}:V(\lambda
+k\kappa )^{GL_{n}}]
\end{equation*}
holds for any integer $k\geq \frac{\left| \nu \right| -\left| \lambda
\right| }{2}=\frac{\left| \mu \right| -\left| \lambda \right| }{2}.$ Since $%
\lambda $ and $\nu $ are partitions, we can write $[V(\nu +k\kappa )^{%
\widehat{G}}:V(\lambda +k\kappa )^{GL_{n}}]=[V(\nu )^{\widehat{G}}:V(\lambda
)^{GL_{n}}]$ (see Remark after Proposition \ref{prop_mul_sum}), this yields
to the desired equality.
\end{proof}

\bigskip

\noindent \textbf{Remark}$\mathrm{:}$ By setting $q=1$ in the above
identity, we obtain the identity 
\begin{equation*}
\frak{D}_{\mu ^{(1)},...,\mu ^{(r)}}^{\lambda ,G}=\sum_{\nu \in \mathcal{P}%
_{n}}[V(\lambda )^{\widehat{G}}:V(\nu )^{GL_{n}}]c_{\mu ^{(1)},...,\mu
^{(r)}}^{\nu }
\end{equation*}
which can also be deduced from the decompositions 
\begin{equation*}
\frak{V}^{G}(\mu ^{(1)})\otimes \cdot \cdot \cdot \otimes \frak{V}^{G}(\mu
^{(r)})\simeq \bigoplus_{\nu \in \mathcal{P}_{n}}V^{GL_{2n}}(\nu )^{\oplus
c_{\mu ^{(1)},...,\mu ^{(r)}}^{\nu }}
\end{equation*}
and 
\begin{equation*}
V^{GL_{2n}}(\nu )\simeq \bigoplus_{\lambda \in \mathcal{P}_{n}}[V(\nu
)^{GL_{2n}}:V(\lambda )^{G}]
\end{equation*}
since we have the duality $[V(\nu )^{GL_{2n}}:V(\lambda )^{G}]=[V(\lambda )^{%
\widehat{G}}:V(\nu )^{GL_{n}}]$ (see Theorem $A_{1}$ of \cite{K}).

\bigskip

\noindent When the $\mu ^{(k)}$'s are rectangular partitions of decreasing
heights, $\frak{D}_{\mu ^{(1)},...,\mu ^{(r)}}^{\lambda ,G}(q)$ has
nonnegative coefficients by (\ref{dualA}) and Theorem \ref{th-shimo}. From
the above theorem and (\ref{express}) we obtain 
\begin{equation}
\frak{D}_{\mu ^{(1)},...,\mu ^{(r)}}^{\lambda ,G}(q)=q^{\tfrac{\left| \mu
\right| -\left| \lambda \right| }{2}}\sum_{\nu \in \mathcal{P}%
_{n}}\sum_{\gamma \in \mathcal{P}_{n}^{_{\diamondsuit _{G}}}}c_{\gamma
,\lambda }^{\nu }K_{\nu ,\mu }^{GL_{n},I}(q)  \label{decf}
\end{equation}
where 
\begin{equation*}
\diamondsuit _{G}=\left\{ 
\begin{tabular}{l}
$(1,1)$ if $G=SO_{2n}$ \\ 
$(2)$ if $G=Sp_{2n}$ \\ 
$(1)$ if $G=SO_{2n+1}$%
\end{tabular}
\right. .
\end{equation*}
For completion, set $\frak{D}_{\mu ^{(1)},...,\mu ^{(r)}}^{\nu
,GL_{n}}(q)=K_{\nu ,\mu }^{GL_{n},I}(q)$ and $\diamondsuit
_{GL_{n}}=\emptyset .\;$Then the four families of polynomials $\frak{D}_{\mu
^{(1)},...,\mu ^{(r)}}^{\lambda ,G}(q)$ can be classified by using the same
symbols $\diamondsuit _{G}=\emptyset ,(1),(1,1),(2)$ labelling the four
families of stable one-dimensional sums.\ Write $X_{\mu ^{(1)},...,\mu
^{(r)}}^{\lambda ,\diamondsuit _{G}}(q)$ for the one-dimensional sum defined
from the affine crystal of kind $\diamondsuit _{G}$ associated to $\mu
^{(1)},...,\mu ^{(r)}$ and $\lambda $.\ By Theorem \ref{th-shimo} we know
that $X_{\mu ^{(1)},...,\mu ^{(r)}}^{\nu ,\emptyset }(q)=K_{\nu ,\mu
}^{GL_{n},I}(q).\;$Then we derive from (\ref{decf}) and Theorem \ref
{th_qdual} that Conjecture $5$ of \cite{sh} giving the decomposition of the
one-dimensional sums of kind $\diamondsuit _{G}=(1),(1,1),(2)$ in terms of
those of kind $\emptyset $ can be reformulated as follows:

\begin{conjecture}
\label{conj_last}The $q$-analogue $\frak{D}_{\mu ^{(1)},...,\mu
^{(r)}}^{\lambda ,G}(q)$ coincide with the stable one-dimensional sum $%
X_{\mu ^{(1)},...,\mu ^{(r)}}^{\lambda ,\diamondsuit _{G}}(q)$ up to the
multiplication by a power of $q:$%
\begin{equation*}
\frak{D}_{\mu ^{(1)},...,\mu ^{(r)}}^{\lambda ,G}(q)=\widetilde{K}_{\widehat{%
\lambda },\widehat{\mu }}^{\widehat{G},I}(q)=q^{\ast }X_{\mu ^{(1)},...,\mu
^{(r)}}^{\lambda ,\diamondsuit _{G}}(q)
\end{equation*}
where $q^{\ast }$ is a power of $q$ depending on the normalization of the
energy (or co-energy) function chosen on this crystal.
\end{conjecture}

\noindent \textbf{Remark}$\mathrm{:}$ This conjecture has been proved in 
\cite{LS} when $\mu ^{(1)},...,\mu ^{(r)}$ are row partitions of decreasing
lengths. In this case the stable one-dimensional sums are stable limits of
Lusztig $q$-analogues.\ Accordingly to the previous conjecture, the stable
one-dimensional sums associated to nonexceptional affine Lie algebras are
stable limits of generalized Lusztig $q$-analogues.

\section{Question}

In Theorem \ref{th_qdual}, the partitions $\mu ^{(1)},...,\mu ^{(r)}$ are
not supposed rectangular.\ Conjecture \ref{conj} suggests that it should
exists, for all classical Lie groups, nonnegative $q$-analogues $\frak{D}%
_{\mu ^{(1)},...,\mu ^{(r)}}^{\lambda ,G}(q)$ of tensor product coefficients
defined from the $r$-tuples $\mathbf{\mu }=(\mu ^{(1)},...,\mu ^{(r)})\;$%
such that $\mu $ is a partition.\ It would be interesting to find a
combinatorial interpretation of the $\frak{D}_{\mu ^{(1)},...,\mu
^{(r)}}^{\lambda ,G}(q)$ when $\mu ^{(1)},...,\mu ^{(r)}$ are not
rectangular.

\end{document}